\documentclass[11pt]{amsart}
\usepackage{color}

\usepackage[top=26mm, bottom=26mm, left=30mm, right=30mm]{geometry}
\usepackage{bm}

\newcommand{\mtx}[1]{\bm{\mathsf{#1}}}

\usepackage{verbatim}
\usepackage{amsthm}
\usepackage{amstext}
\usepackage{graphicx}
\usepackage{amssymb}

\makeatletter

\begin{document}

\title[High-order approximation of evolution operators]{A high-order scheme for solving wave propagation problems via
the direct construction of an approximate time-evolution operator}

\author{T. S. Haut}
\author{T. Babb}
\author{P. G. Martinsson}
\author{B. A. Wingate}

\maketitle

\begin{center}
\begin{minipage}{0.85\textwidth}\small
\textbf{Abstract} The manuscript presents a technique for efficiently solving the classical
wave equation, the shallow water equations, and, more generally,
equations of the form $\partial u /\partial t = \mathcal{L}u$, where
$\mathcal{L}$ is a skew-Hermitian differential operator. The idea is to explicitly
construct an approximation to the time-evolution operator
$\exp(\tau\mathcal{L})$ for a relatively large time-step $\tau$.
Recently developed techniques for approximating oscillatory scalar functions by
rational functions, and accelerated algorithms for computing functions of
discretized differential operators are exploited. Principal advantages of
the proposed method include: stability even for large time-steps, the
possibility to parallelize in time over many characteristic wavelengths, and large speed-ups over existing methods
in situations where simulation over long times are required

\vspace{1mm}

Numerical examples involving the 2D rotating shallow water equations and the
2D wave equation in an inhomogenous medium are presented, and the method
is compared to the $4$th order Runge-Kutta (RK4) method and to the use of
Chebyshev polynomials. The new method achieved high accuracy over long
time intervals, and with speeds that are orders of magnitude faster
than both RK4 and the use of Chebyshev polynomials.
\end{minipage}
\end{center}

\section{Introduction}

\subsection{Problem formulation}

We present a technique for solving a class of linear hyperbolic problems
\begin{equation}
\label{eq:basic}
\left\{\begin{aligned}
\frac{\partial \mathbf{u}}{\partial t}(\mathbf{x},t) =&\  \mathcal{L}\mathbf{u}(\mathbf{x},t),\qquad&&\mathbf{x} \in \Omega,\quad t > 0,\\
\mathbf{u}(\mathbf{x}, 0) =&\ \mathbf{u}_{0} \left( \mathbf{x} \right)&&\mathbf{x} \in \Omega .
\end{aligned}\right.
\end{equation}
Here $\mathbf{u}$ is a possibly vector valued function and  $\mathcal{L}$ is a skew-Hermitian differential operator (see the end of this Section for the method's scope).
The technique is demonstrated on the 2D rotating shallow water equations, as well as the variable coefficient wave equation.

The basic approach is classical, and involves the construction of a rational approximation to the
time evolution operator $\exp(\tau \mathcal{L})$ in the form
$$
\exp(\tau \mathcal{L}) \approx \sum_{m=-M}^{M}b_{m}\,\bigl(\tau \mathcal{L} - \alpha_{m}\bigr)^{-1},
$$
where the time-step $\tau$ is fixed in advance and $M$ scales linearly in $\tau$. Once the time-step $\tau$ has been fixed,
an approximate solution at times $\tau,\,2\tau,\,3\tau,\dots$ can be obtained via repeated application of
the approximate time-stepping operator, since $\exp(n \tau \mathcal{L})  = \bigl(\exp(\tau\mathcal{L})\bigr)^{n}$.
The computational profile of the method is that it takes a moderate amount of work to construct the initial
approximation to $\exp(\tau\mathcal{L})$, but once it has been built, it can be applied very rapidly, even
for large $\tau$.

The efficiency of the proposed scheme is enabled by (i) a novel method  \cite{D-B-H-M:2013} for constructing near
optimal rational approximations
to oscillatory functions such as $e^{i x}$ over arbitrarily long intervals, and by (ii) the development \cite{MARTIN:2013} of a
high-order accurate and stable method for pre-computing approximations to operators
of the form $(\tau \mathcal{L} - \alpha_{m})^{-1}$. The near optimality of the rational approximations ensures that
the number $2 M +1$ of terms needed for a given accuracy is typically much smaller than standard
methods that rely on polynomial or rational approximations of $\mathcal{L}$.

The proposed scheme has several advantages over typical methods, including the absence of stability constraints on
the time step $\tau$ in relation to the spatial discretization,
the ability to parallelize in time over many characteristic wavelengths (in addition to any spatial parallelization),
and great acceleration when integrating equation (\ref{eq:basic}) for long times or for multiple initial conditions
(e.g.~when employing an exponential integrator on a nonlinear evolution equation, cf.~Section \ref{sec:generalizations}).
A drawback of the scheme is that it is more memory intensive than standard techniques.

We restrict the scope of this paper to when the application $(\tau \mathcal{L} - \alpha_m)^{-1} \mathbf{u}_0$ can be reduced to the solution of an  elliptic-type PDE for one of the unknown variables. This situation arises in geophysical fluid applications (among others), including  the rotating primitive equations that are that are used for climate simulations. In this context, the ability to efficiently solve (\ref{eq:basic}) can be used to construct efficient schemes for the fully nonlinear evolution equations in the presence of time scale separation (see \cite{HAU-WIN:2013}). However, the direct solver presented in Section~\ref{sec:Spectral-element-discretization} is quite general, and in principle can be extended to first order linear systems of hyperbolic PDEs with little modification (though such an extension is speculative and, in particular, has not been tried).


\subsection{Time discretization}
\label{sec:intro_time}

In order to time-discretize (\ref{eq:basic}), we fix a time-step $\tau$ (the choice of which is discussed shortly),
a requested precision $\delta > 0$,
and ``band-width'' $\Lambda \in (0,\infty)$ which specifies the spatial resolution
(in effect, the scheme will accurately capture eigenmodes of $\mathcal{L}$ whose
eigenvalues $\lambda$ satisfy $|\lambda| \leq \Lambda$). We then use an improved
version of the scheme of \cite{D-B-H-M:2013} to construct a rational
function,
\begin{equation}
\label{eq:def_RM}
R_{M}(ix) = \sum_{m=-M}^{M} \frac{b_{m}}{ \bigl(ix - \alpha_{m} \bigr) },
\end{equation}
such that
\begin{equation}
\label{eq:RM_prop1}
\bigl| e^{i x} - R_{M}(ix)\bigr| \leq \delta,\qquad x \in [-\tau \Lambda, \tau \Lambda],
\end{equation}
and
\begin{equation}
\label{eq:RM_prop2}
|R_{M}(ix)| \leq 1,\qquad x \in \mathbb{R}.
\end{equation}
It now follows from
(\ref{eq:RM_prop1}) and (\ref{eq:RM_prop2}) that if we approximate $\exp(t\mathcal{L})$
by $R_{M}(\tau \mathcal{L})$, the approximation error satisfies
\begin{equation}
\label{eq:error in exp(t*L)}
\left\Vert e^{\tau \mathcal{L}} \mathbf{u}_{0}-
\sum_{m=-M}^{M}b_{m}\,\bigl(\tau \mathcal{L} - \alpha_{m}\bigr)^{-1} \mathbf{u}_{0}\right\Vert
\leq\delta\left\Vert \mathbf{u}_{0}\right\Vert +
2\left\Vert \mathbf{u}_{0}-\mathcal{P}_{\Lambda} \mathbf{u}_{0}\right\Vert,
\end{equation}
where $\mathcal{P}_{\Lambda}$ projects functions onto the subspace spanned
by eigenvectors of $\mathcal{L}$ with modulus at most $\Lambda$. Here the
only property of $\mathcal{L}$ that we use is that $\mathcal{L}$ is
skew-Hermitian, and hence has a complete spectral decomposition with a
purely imaginary spectrum.

The bound (\ref{eq:RM_prop2}) ensures
that the repeated application of $R_{M}(\tau \mathcal{L})$  is stable on the entire imaginary axis. It also turns out that
the number $2 M+1$ of terms needed in the rational approximation in (\ref{eq:RM_prop1}) is close to optimally small (for the given
accuracy $\delta$).

The scheme described above allows a great deal of freedom in the choice of the
time step $\tau$. While classical methods typically require
the time step to be a small fraction of the characteristic wavelength,
we have freedom to let $\tau$ cover a large number of characteristic wavelengths. Therefore,
the scheme is well suited to parallelization in time, since
all the inverse operators in the approximation of the operator exponential can be applied independently.
In fact, the only constraint on the size of $\tau$ is on
the memory available to store the representations of the inverse
operators (as explained in Section \ref{subsec: intro, pre-comp of inverses}, the memory required
for each inverse scales linearly in the number of spatial discretization parameters, up to a logarithmic factor).

\subsection{Pre-computation of rational functions of $\mathcal{L}$} \label{subsec: intro, pre-comp of inverses}

The time discretization technique described in Section \ref{sec:intro_time}
requires us to build explicit approximations to differential
operators on the domain $\Omega$ such as  $(\tau \mathcal{L} - \alpha_{m})^{-1}$.
We do this using a variation of the technique described in \cite{MARTIN:2013}.
A variety of different domains can be handled, but for simplicity, suppose that
$\Omega$ is a rectangle. The idea is to tessellate $\Omega$ into a collection
of smaller rectangles, and to put down a tensor product grid of Chebyshev nodes
on each rectangle, as shown in Figure~\ref{fig:grid}. A function is represented via tabulation
on the nodes, and then $\mathcal{L}$ is discretized via standard spectral collocation
techniques on each patch. The patches are glued together by enforcing continuity
of both function values and normal derivatives. This discretization results in a block sparse coefficient
matrix, which can rapidly be inverted via a procedure very similar to the classical
nested dissection technique of George \cite{george_1973}. The resulting inverse is
dense but ``data-sparse,'' which is to say that it has internal structure that allows
us to store and apply it efficiently.

In order to describe the computational cost of the direct solver, let
$N$ denote the number of nodes in the spatial discretization. For a problem in
two dimensions, the ``build stage'' of the proposed scheme constructs $2M+1$
data-sparse matrices $\{\mtx{A}_{m}\}_{m=-M}^{M}$ of size $N\times N$, where
each $\mtx{A}_{m}$ approximates $(\tau \mathcal{L} - \alpha_{m})^{-1}$.
The build stage has asymptotic cost $\mathcal{O}(M\,N^{1.5})$, and storing the matrices
requires $O(M\,N\,\log(N))$ memory. The cost of applying a matrix $\mtx{A}_{m}$ is
$\mathcal{O}(N\,\log\left(N\right))$. (We remark that the cost of building the matrices $\{\mtx{A}_{m}\}_{m=-M}^{M}$ can often be
accelerated to optimal $\mathcal{O}(M\,N)$ complexity \cite{2013_martinsson_DtN_linearcomplexity},
but since the pre-factor in the $\mathcal{O}(M\,N^{1.5})$ bound is quite small,
such acceleration would have negligible benefit for the problem sizes
under consideration here.) Section \ref{sec:Spectral-element-discretization}
describes the inversion procedure in more detail.

We remark that the spatial discretization procedure we use does not explicitly enforce
that the discrete operator is exactly skew-Hermitian. However, the fact that the spatial
discretization is done to very high accuracy means that it is in practice very nearly so.
Numerical experiments indicate that the scheme as a whole is stable in every regime
where it was tested. 




\subsection{Comparison to existing approaches}

The approach of using proper rational approximations for applying
matrix exponentials has a long history. In the context of operators with negative spectrum (e.g. for parabolic-type PDEs),
many authors have discussed how to compute efficient
rational approximations to the decaying exponential $e^{-x}$, including
using Cauchy's integral formula coupled with Talbot quadrature (cf.~\cite{TR-WE-SC:2006}), and optimal rational approximations via the
Carath\'eodory-Fejer method (cf.~\cite{TR-WE-SC:2006}) or the Remez algorithm \cite{CO-ME-VA:1969}. However, such methods are less effective
(or not applicable) when applied to approximating oscillatory functions
such as $e^{ix}$ over long intervals. For computing functions of parabolic-type linear operators,
the approach of combining rational approximations and compressed representations
of the solution operators using so-called $\mathcal{H}$-matrices has been proposed
in \cite{GA-HA-KH:2005}.

Common approaches for applying the exponential of skew-Hermitian operators
include high-order time-stepping methods, scaling-and-squaring coupled with  Pad\'e approximations (cf. \cite{HIGHAM:2005}) or
Chebyshev polynomials (cf. \cite{BER-VIA:2000}), and polynomial or
rational Krylov methods (cf. \cite{HOC-LUB:1997} and \cite{GUTTEL:2012}). All these methods
iteratively build up rational or polynomial approximations to the operator exponential, and correspondingly
approximate the spectrum $e^{i\omega_{n}\tau}$
of $e^{\tau \mathcal{L}}$ with polynomials or rationals. Therefore, the near optimality of (\ref{eq:def_RM})
and the speed of applying the inverse operators in (\ref{eq:error in exp(t*L)})
will generally translate into high efficiency relative to standard methods. In contrast to these standard approaches,
the method proposed in this paper can also be trivially parallelized in time over many characteristic wavelengths.

In addition to approaches that rely on polynomial or rational
approximations, let us mention two alternative approaches for time-stepping on wave propagation problems.
The authors in \cite{BEY-SAN:2005} combine separated representations of multi-dimensional operators, partitioned low
rank compressions of matrices, and (near) optimal quadrature nodes for band-limited functions, in order to compute compressed representations
of the operator exponential over $1-2$ characteristic wavelengths. Along different lines, the authors in \cite{DEM-YIN:2009}
use wave atoms to construct compressed representations of the (short
time) operator exponential, and in particular can bypass the CFL constraint.


\subsection{Outline of manuscript}

The paper is organized as follows.
In Section~\ref{sec:Spectral-element-discretization}, we briefly
describe the direct solver in \cite{MARTIN:2013}. We then discuss in Section~\ref{sec:constructing-rational-approximations}
a technique for constructing efficient rational approximations
of general functions, and specialize to the case of approximating
the exponential $e^{i x}$ and the phi-functions for exponential integrators \cite{COX-MAT:2002}.
In Section~\ref{sec:Examples}, we present applications of
the method for both the 2D rotating shallow water equations and the
2D wave equation in inhomogenous medium. In particular, we compare the accuracy and
efficiency of this approach against $4$th order Runge-Kutta and the
Chebyshev polynomial method (in our comparisons, we use the same spectral
element discretization). Finally, Appendix~\ref{Appendix: error bounds} contains
error bounds for the rational approximations constructed here.

\section{Spectral element discretization\label{sec:Spectral-element-discretization}}

This section describes how to efficiently compute a highly accurate approximation to
the inverse operator $\left( \mathcal{L} - \alpha \right)^{-1}$, where $\mathcal{L}$
is a skew-Hermitian operator. As mentioned in the introduction, we restrict our
discussion to environments where application of the inverse can be reformulated as
a scalar elliptic problem. This reformulation procedure is illustrated for the classical
wave equation and for the shallow water equations in Section \ref{sec:reformulation}.
Section \ref{sec:discretization} describes a high-order multidomain spectral discretization
procedure for the elliptic equation. Section \ref{sec:directsolver} describes a direct
solver for the system of linear equations arising upon discretization.

\subsection{Reformulation as an elliptic problem}
\label{sec:reformulation}
In many situations of practical interest, the task of solving a hyperbolic equation
$(\mathcal{L} - \alpha)u = f$, where $\mathcal{L}$ is a skew-Hermitian operator,
can be reformulated as an associated elliptic problem. In this section, we illustrate
the idea via two representative examples. Example $1$ is of particular relevance to geophysical fluid
applications, which serve as a major motivation of this algorithm.

\vspace{2mm}

\noindent
\textit{Example 1 --- the shallow water equation:}
We consider the rotating shallow water equations,
\begin{eqnarray}
\mathbf{v}_{t} & = & -fJ\mathbf{v}+\nabla\eta,\label{eq:RSW equations}\\
\eta_{t} & = & \nabla\cdot\mathbf{v},\nonumber
\end{eqnarray}
where $\mathbf{v}\left(\mathbf{x}\right)=\left(v_{1}\left(\mathbf{x}\right),v_{2}\left(\mathbf{x}\right)\right)$
denotes the fluid velocity, $\eta\left(\mathbf{x}\right)$ denotes
perturbed surface elevation, $f$ is the (possibly spatially varying)
Coriolis frequency, and
\[
J=\left(\begin{array}{cc}
0 & 1\\
-1 & 0\end{array}\right).
\]
On the sphere, $f=2\Omega\sin\phi$; on the plane, $f$ is constant.
We write system (\ref{eq:RSW equations}) in the form
\[
\mathbf{u}_{t}=\mathcal{L}\mathbf{u},
\]
where
\begin{equation}
\mathcal{L}\left(\begin{array}{c}
\mathbf{v}\\
\eta\end{array}\right)=\left(\begin{array}{c}
-fJ\mathbf{v}+\nabla\eta\\
\nabla\cdot\mathbf{v}\end{array}\right).\label{eq:forward oper, RSW}
\end{equation}
Although we only consider the case when the Coriolis frequency is
constant, the method generalizes to non-constant coefficient $f$
(see also the example in the next section) and is of particular relevance
for a spectral element discretization on the cubed sphere.

In order to apply the method in this paper, we use the standard fact
(cf. \cite{PAL-SIG:2011}) that if
\begin{equation}
\left(\mathcal{L}-\alpha\right)\left(\begin{array}{c}
\mathbf{v}\\
\eta \end{array}\right)=\left(\begin{array}{c}
\mathbf{v_0}\\
\eta_0\end{array}\right),\label{eq:inverse oper, RSW eqns}
\end{equation}
then $\eta$ satisfies the elliptic equation
\begin{equation}
\nabla\cdot\left(\mathcal{A}_{\alpha}\nabla\eta\right)-\alpha \eta=\eta_0+H\nabla\cdot\mathcal{A}_{\alpha}\mathbf{v_0}.\label{eq:elliptic solve for RSW eqns}\end{equation}
\mbox{} Here $\mathcal{A}_{\alpha}$ is defined by \[
\mathcal{A}_{\alpha}=\frac{1}{\alpha^{2}+f^{2}}\left(\begin{array}{cc}
\alpha & f\\
-f & \alpha\end{array}\right).\]
Once $\eta$ is computed, $\mathbf{v}$ can be obtained directly,\begin{equation}
\mathbf{v}=-\mathcal{A}_{\alpha}\mathbf{v_0}+\mathcal{A}_{\alpha}\nabla\eta.\label{eq:v_m, RSW}\end{equation}
When $f$ is constant, equation (\ref{eq:elliptic solve for RSW eqns})
reduces to
\begin{equation}
\left(\Delta-\frac{\alpha^{2}+f^{2}}{c^{2}}\right)\eta=\frac{\alpha^{2}+f^{2}}{c^{2}\alpha}\left(\eta_0+H\nabla\cdot\left(\mathcal{A}_{\alpha}\mathbf{v_0}\right)\right).
\label{eq:elliptic solve, RSW}
\end{equation}

\vspace{2mm}

\noindent
\textit{Example 2 --- the wave equation:}
Consider the wave propagation problem
\begin{equation}
\label{eq:wave equation}
u_{tt}=\kappa\Delta u,\,\,\,\,\mathbf{x}\in\left[0,1\right]\times\left[0,1\right],
\end{equation}
where $\kappa\left(\mathbf{x}\right)>0$ is a smooth function, the
initial conditions $u\left(\mathbf{x},0\right)$ and $u_{t}\left(\mathbf{x},0\right)$
are prescribed, and periodic boundary conditions are used.

In order to apply the method in this paper, we reformulate (\ref{eq:wave equation})
as a first order system in both time and space by defining $v=u_{t}$,
$w=u_{x}$, and $z=u_{y}$. Then we have that\begin{equation}
\left(\begin{array}{c}
w_{t}\\
z_{t}\\
v_{t}\end{array}\right)=\left(\begin{array}{ccc}
0 & 0 & \partial_{x}\\
0 & 0 & \partial_{y}\\
\kappa\partial_{x} & \kappa\partial_{y} & 0\end{array}\right)\left(\begin{array}{c}
w\\
z\\
v\end{array}\right),
\label{eq:forward oper, inhom}
\end{equation}
with initial conditions
\[
v\left(\mathbf{x},0\right)=u_{0}\left(\mathbf{x}\right),\,\,\, w\left(\mathbf{x},0\right)=\frac{\partial u_{0}}{\partial x}\left(\mathbf{x}\right),\,\,\, z\left(\mathbf{x},0\right)=\frac{\partial u_{0}}{\partial y}\left(\mathbf{x}\right).
\]
Here the scalar function $u$ to the original system (\ref{eq:wave equation})
can be recovered after the final time step by solving the elliptic
equation $\Delta u=w_{x}+z_{y}$.

To apply the method in this paper, we compute the solution to
\begin{equation}
\label{eq:operator for inhomog, inverse}
\left(\mathcal{L}-\alpha\right)\left(\begin{array}{c}
w\\
z\\
v\end{array}\right)=\left(\begin{array}{c}
v_{x}-\alpha w\\
v_{y}-\alpha z\\
\kappa\left(w_{x}+z_{y}\right)-\alpha v\end{array}\right)=\left(\begin{array}{c}
w_{0}\\
z_{0}\\
v_{0}\end{array}\right)\end{equation}
as follows. First, solving for $w$ and $z$ in terms of $v$,
\begin{equation}
\label{eq:w and z eqns}
w=\frac{1}{\alpha}\left(v_{x}-w_{0}\right),\,\,\,\, z=\frac{1}{\alpha}\left(v_{y}-z_{0}\right),
\end{equation}
it is straightfoward to show that
\begin{equation}
\label{eq:v for inhomog, inverse}
\left(\Delta-\alpha^{2}\kappa^{-1}\right)v=\alpha\kappa^{-1}v_{0}+\frac{\partial w_{0}}{\partial x}+\frac{\partial z_{0}}{\partial y}.
\end{equation}
Once $v$ is known, $w$ and $z$ can then be computed directly via
(\ref{eq:w and z eqns}).

\subsection{Discretization}
\label{sec:discretization}
In this section, we describe a high-order accurate discretization
scheme for elliptic boundary value problems such as (\ref{eq:elliptic solve, RSW})
and (\ref{eq:v for inhomog, inverse}) which arise in the solution of
hyperbolic evolution equations. Specifically, we describe the
solver for a boundary value problem (BVP) of the form
\begin{equation}
\label{eq:elliptic}
\mathcal{B}u\left(\mathbf{x}\right)=f\left(\mathbf{x}\right),\,\,\,
\mathbf{x}\in\Omega,
\end{equation}
where $\mathcal{B}$ is an elliptic differential operator. To keep
things simple, we consider only square domains $\Omega = [0,1]^{2}$,
but the solver can easily be generalized to other domains. The
solver we use is described in detail in \cite{2013_martinsson_DtN_tutorial},
our aim here is merely to give a high-level conceptual description.

\begin{figure}
\includegraphics[width=115mm]{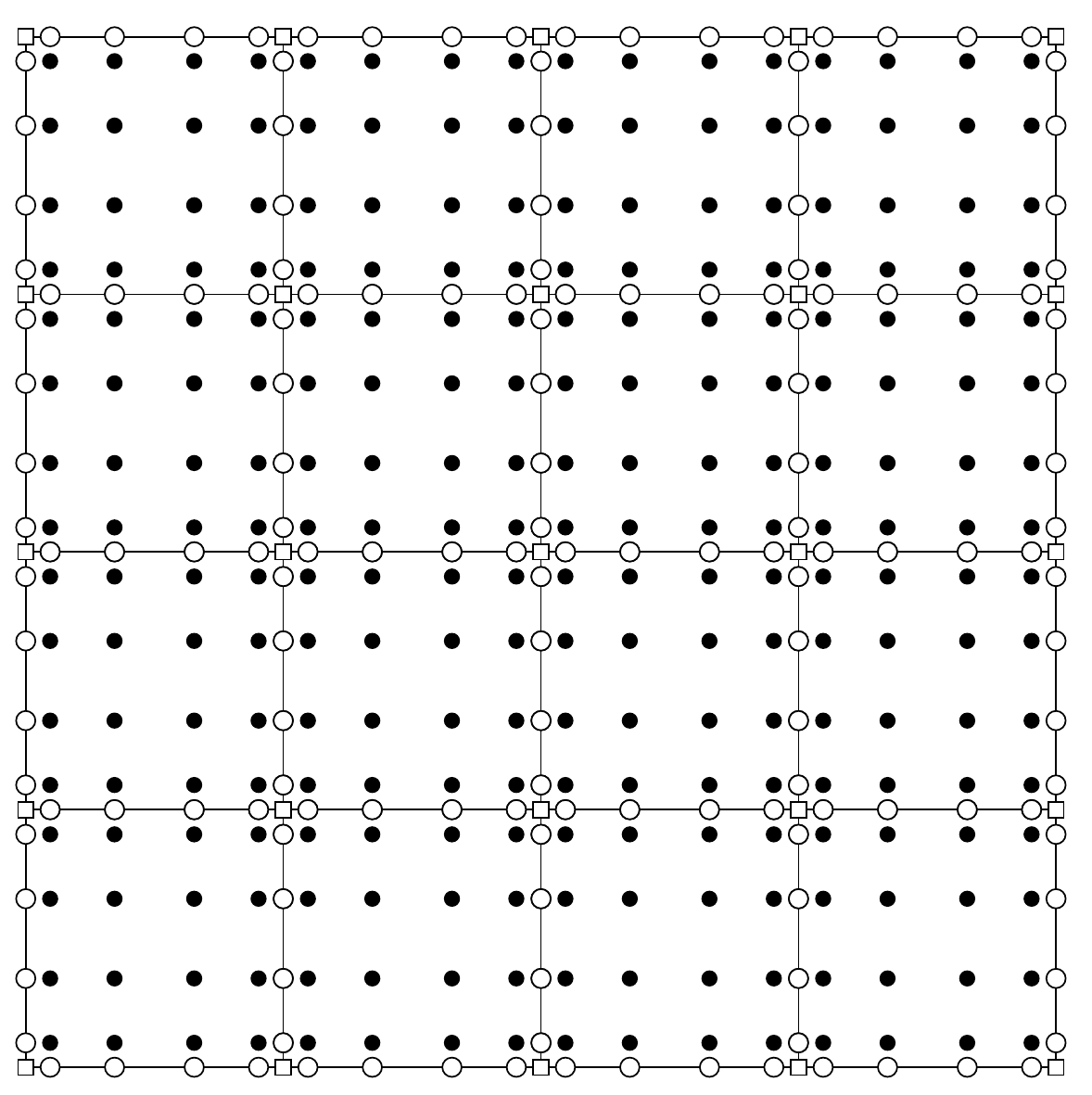}
\caption{Illustration of the grid of points $\{\mathbf{x}_{j}\}_{j=1}^{N}$ introduced
to discretize (\ref{eq:elliptic}) in Section \ref{sec:discretization}. The figure
shows a simplified case involving $4\times 4$ squares, each holding a $6\times 6$
local tensor product grid of Chebyshev nodes. The PDE (\ref{eq:elliptic}) is
enforced via collocation using spectral differentiation on each small square
at all solid (``internal'') nodes. At the hollow (``boundary'') nodes, continuity
of normal fluxes is enforced.}
\label{fig:grid}
\end{figure}

The PDE (\ref{eq:elliptic}) is discretized using a multidomain
spectral collocation method. Specifically, we split the square
$\Omega$ into a large number of smaller squares (or rectangles),
and then put down a tensor product grid of $p \times p$
Chebyshev nodes on each small square, see Figure \ref{fig:grid}.
The parameter $p$ is chosen so that dense computations involving
matrices of size $p^{2}\times p^{2}$ are cheap ($p=20$ is often a good choice).
Let $\{\mathbf{x}_{j}\}_{j=1}^{N}$ denote the total set of nodes.
Our approximation to the solution $u$ of (\ref{eq:elliptic}) is
then represented by a vector $\mathbf{u} \in \mathbb{C}^{N}$,
where the $j$'th entry is simply an approximation to the function
value at node $\mathbf{x}_{j}$, so that $\mathbf{u}(j) \approx u(\mathbf{x}_{j})$.
The discrete approximation to (\ref{eq:elliptic}) then takes the form
\begin{equation}
\label{eq:Bu=f}
\mathbf{B}\mathbf{u} = \mathbf{f},
\end{equation}
where $\mathbf{B}$ is an $N\times N$ matrix. The $j$'th row of
(\ref{eq:Bu=f}) is associated with a collocation condition for
node $\mathbf{x}_{j}$. For all $j$ for which $\mathbf{x}_{j}$
is a node in the \textit{interior} of a small square (filled
circles in Figure \ref{fig:grid}), we directly
enforce (\ref{eq:elliptic}) by replacing all differentiation
operators by spectral differentiation operators on the local
$p\times p$ tensor product grid. For all $j$ for which $\mathbf{x}_{j}$
lies on a \textit{boundary} between two squares (hollow squares
in Figure \ref{fig:grid}), we enforce that normal fluxes across
the boundary are continuous, where the fluxes from each side of
the boundary are evaluated via spectral differentiation on the
two patches (corner nodes need special treatment, see \cite{2013_martinsson_DtN_tutorial}).

\subsection{Direct solver}
\label{sec:directsolver}
The discrete linear system (\ref{eq:Bu=f}) arising from discretization
of (\ref{eq:elliptic}) is block-sparse. Since it has the typical
sparsity pattern of a matrix discretizing a 2D differential operator,
it is possible to compute its LU factorization in $O(N^{1.5})$ operations
using a nested dissection ordering of the nodes \cite{DU-ER-RE:1986,george_1973}
that minimizes fill-in.
Once the LU-factors have been computed, the cost of a linear solve is
$O(N \log N)$. In the numerical computations presented in Section \ref{sec:Examples},
we use a slight variation of the nested-dissection algorithm that was introduced
in \cite{MARTIN:2013} for the case of homogeneous equations. The extension to the
situation involving body loads is straight-forward, see \cite{2013_martinsson_DtN_tutorial}.

We note that by exploiting internal structure in the dense sub-matrices that
appear in the factors of $\mathbf{B}$ as the factorization proceeds, the complexity
of both the factorization and the solve stages can often be reduced to optimal $O(N)$
complexity \cite{2013_martinsson_DtN_linearcomplexity}. However, for the problem sizes 
considered in this manuscript, there would be little practical gain to
implementing this more complex algorithm.

\section{constructing rational approximations\label{sec:constructing-rational-approximations}}

We now discuss how to construct efficient rational approximations
to general smooth functions $f\left(x\right)$. For concreteness,
we consider approximating the phi functions \[
\varphi_{0}\left(x\right)=e^{ix},\,\,\,\varphi_{1}\left(x\right)=\frac{e^{ix}-1}{ix},\,\,\,\varphi_{2}\left(x\right)=\frac{e^{ix}-ix-1}{ix^{2}},\]
that arise for high-order exponential integrators (cf. \cite{TR-WE-SC:2006}).
By considering the real and imaginary components separately, we assume
that $f\left(x\right)$ is real-valued (it turns out that the poles
in the approximation will be the same for the real and imaginary components).

The construction proceeds in two steps; the second step is actually
a pre-computation and need only be done once, but is presented last
for clarity. First, we construct an approximation to $f\left(x\right)$
by sums of shifted Gaussians $\psi_{h}\left(x\right)=\left(4\pi\right)^{-1/2}e^{-x^{2}/\left(4h^{2}\right)}$
(see Section~\ref{sub:Gaussian-approximations-to} for details),
\begin{equation}
\left|f\left(x\right)-\sum_{-M}^{M}b_{m}\psi_{h}\left(x+nh\right)\right|\leq\delta_{1},\,\,\,\,-\Lambda\leq x\leq\Lambda.\label{eq:intro, approx of f(x) via Gaussians}\end{equation}
Here $h$ is inversely proportional to the bandlimit of $f\left(x\right)$,
and $M$ controls the interval $\Lambda$ over which the approximation
is valid (roughly $\left|x\right|\lesssim Mh$). When $f\left(x\right)=e^{ix}$,
the coefficients are explicitly given by $c_{m}=\left(\widehat{\psi_{h}}\left(1\right)/h\right)e^{-2\pi inh}$,
and the approximation is remarkably accurate (see \ref{eq:error bound for e^ix, partial-1}
for error bounds). Second, using the approach in \cite{D-B-H-M:2013},
a rational approximation to $\psi_{1}\left(x\right)=\left(4\pi\right)^{-1/2}e^{-x^{2}/4}$
is constructed over the real line (see Section~\ref{sub:Rational-approximation-to}
for details),\begin{equation}
\left|\psi_{1}\left(x\right)-2\text{Re}\left(\sum_{j=-L}^{L}\frac{a_{j}}{ix-\left(\mu+ij\right)}\right)\right|\leq\delta_{2},\,\,\,\, x\in\mathbb{R}.\label{eq:intro, rational approx to Gaussian}\end{equation}
Notice that the imaginary parts of the poles in the above approximation
are integer multiples $j=0,\pm1,\ldots,\pm L$. For $L=11$, we construct
$\mu$ and coefficients $a_{j}$ such that the $L^{\infty}$ approximation
error $\delta_{2}$ satisfies $\delta_{2}<10^{-12}$ (see Table~\ref{tab:Coefficients for rational approx Gaussian}).
Finally, combining (\ref{eq:intro, approx of f(x) via Gaussians})
and (\ref{eq:intro, rational approx to Gaussian}), we obtain a rational
approximation to $f\left(x\right)$,\[
\left|f\left(x\right)-2\text{Re}\left(\sum_{n=-M-L}^{M+L}\frac{c_{n}}{ix-h\left(\mu+in\right)}\right)\right|\leq\delta_{1}+2\left(M+L\right)\delta_{2}.\]
Here the coefficients $c_{n}$ are given by \[
c_{n}=h\sum_{k=L_{1}}^{L_{2}}a_{k}b_{n-k},\]
where \[
L_{1}\left(n\right)=\max\left(-L,n-M\right),\,\,\,\, L_{2}\left(n\right)=\max\left(-L,n-M\right).\]

Importantly, constructing the rational approximation (\ref{eq:intro, rational approx to Gaussian})
to $\psi\left(x\right)$ need only be done once. In particular, once
$\mu$ and the coefficients $a_{j}$ are pre-computed, rational approximations
to general functions $f\left(x\right)$ over arbitrarily long
spatial intervals can be obtained with minimal effort, as discussed
in Section~\ref{sub:Gaussian-approximations-to}. We present $\mu$,
and the coefficients $a_{j}$, $j=-11,\ldots,11$, in Table~\ref{tab:Coefficients for rational approx Gaussian},
which are sufficient to yield an $L^{\infty}$ error $\delta_{1}\approx7\times10^{-13}$
in (\ref{eq:intro, rational approx to Gaussian}) .

Using the reduction algorithm in \cite{HAU-BEY:2012}, we
find that the rational approximation constructed for $e^{ix}$ is
close to optimal in the $L^{\infty}$ norm, for a given accuracy $\delta$
and spatial cutoff $A$. In fact, the construction in this paper uses
only $1.2$ times more poles than the near optimal rational approximation
obtained from \cite{HAU-BEY:2012} (when $\delta=10^{-10}$ and $A=56 \pi$,
which we use in our numerical experiments). We note that the residues
corresponding to this near optimal approximation can be very large
and, for this reason, we prefer to use the sub-optimal approximation
instead.

As clarified in Sections~\ref{sub:Gaussian-approximations-to}~and~\ref{sub:Rational-approximation-to},
the same poles can be used to approximate multiple functions with
the same bandlimit. For example, we can use the same poles to approximate
all functions $e^{2\pi itx}$, for $0\leq t\leq1$, since all these
functions have bandlimit less than or equal to $e^{2\pi ix}$; the
dependence on $t$ is only through the coefficients, which are given
explicitly by $c_{m}=\left(\widehat{\psi_{h}}\left(t\right)/h\right)e^{-2\pi inth}$.
In particular, the poles $\alpha_{m}=h\left(\mu+im\right)$ are independent
of $t$ and yield uniformly accurate approximations to $e^{itx}$
on the same interval $\left[-\Lambda,\Lambda\right]$. This observation
enables the efficient computation of multiple operator exponentials
$e^{s_{k}\mathcal{L}}\mathbf{u}_{0}$, for $s_{k}=tk/L$, using the
same computed solutions $\left(t\mathcal{L}-\alpha_{m}\right)^{-1}\mathbf{u}_{0}$,
$m=1,\ldots,M$. A similar comment applies to the phi-functions from
exponential integrators.

Generally, any rational approximation to $e^{ix}$ (or more general
functions) must share the same number of zeros within the interval
of interest; in particular, since the rational approximation can be
expressed as a quotient of polynomials, it is therefore subject to
the Nyquist constraint. However, one advantage of this approximation
method is that it allows efficient rational approximations of functions
that are spatially localized. In fact, since the approximation (\ref{eq:intro, approx of f(x) via Gaussians})
involves highly localized Gaussians, the subsequent rational approximations
are able to represent spatially localized functions as well as highly
oscillatory functions using (perhaps a subset) of the same collection
of poles. This allows the ability
to take advantage of spectral gaps (e.g. from scale separation between
fast and slow waves) and possibly bypass the Nyquist constraint under
certain circumstances.

\subsection{Gaussian approximations to a general function\label{sub:Gaussian-approximations-to}}

We discuss how to construct the approximation (\ref{eq:intro, approx of f(x) via Gaussians}).
To do so, we choose $h$ small enough that the function $\hat{f}\left(\xi\right)$
is zero (or approximately so) outside the interval $\left[-1/\left(2h\right),1/\left(2h\right)\right]$.
Then we can expand $\hat{f}\left(\xi\right)/\widehat{\psi_{h}}\left(\xi\right)$
in a Fourier series, \begin{equation}
\frac{\hat{f}\left(\xi\right)}{\widehat{\psi_{h}}\left(\xi\right)}=\sum_{-\infty}^{\infty}c_{m}e^{2\pi imh\xi},\label{eq:eqn for c_m, fourier domain}\end{equation}
where \[
c_{m}=h\int_{-1/\left(2h\right)}^{1/\left(2h\right)}e^{-2\pi imh\xi}\frac{\hat{f}\left(\xi\right)}{\widehat{\psi_{h}}\left(\xi\right)}d\xi.\]
Transforming (\ref{eq:eqn for c_m, fourier domain}) back to the spatial
domain, we have that \[
f\left(x\right)=\sum_{-\infty}^{\infty}c_{m}\psi_{h}\left(x+mh\right).\]
Notice that the functions $\psi_{h}\left(x+mh\right)$ are tighly
localized in space, and truncating the above series from $-M$ to
$M$ yields accurate approximations for $-\left(M-b\right)hx<x<\left(M-b\right)hx$,
where $b>0$ is a small number that is related to the decay of $\psi_{h}\left(x\right)$.
We remark that the authors in \cite{MAZ-SCH:1996} discuss a related
method of constructing quasi-interpolating representations via sums
of Gaussians (see \cite{MAZ-SCH:2007} for a comprehensive survey).

Specializing to the case when $f\left(x\right)=e^{2\pi ix}$, we have
that $\hat{f}\left(\xi\right)=\delta\left(\xi-1\right)$, and so the
coefficients $c_{m}$ are given by \begin{equation}
c_{m}=\frac{h}{\widehat{\psi_{h}}\left(1\right)}e^{-2\pi imh}.\label{eq:coeff. for exp, gauss}\end{equation}
Similarly, for functions $\varphi_{1}\left(x\right)$ and $\varphi_{2}\left(x\right)$,
the coefficients $c_{m}$ can be obtained numerically using the fact
that \[
\widehat{\phi_{1}}\left(\xi\right)=\begin{cases}
2\pi, & \,\,\,\,-\frac{1}{2\pi}\leq\xi\leq0,\\
0, & \,\,\,\,\text{otherwise.}\end{cases}\]
and\[
\widehat{\phi_{2}}\left(\xi\right)=\begin{cases}
\left(2\pi\right)^{2}\left(\xi+\frac{1}{2\pi}\right), & \,\,\,\,-\frac{1}{2\pi}\leq\xi\leq0,\\
0, & \,\,\,\,\text{otherwise.}\end{cases}\]
For example, the coefficients $c_{m}$ for e.g. $\phi_{1}\left(x\right)$
can be computed via discretization of the integral, \[
c_{m}=h\int_{-1/\left(2\pi\right)}^{0}e^{-2\pi imh\xi}\frac{e^{-2\pi imh\xi}}{\widehat{\psi_{h}}\left(\xi\right)}d\xi.\]

\begin{figure}
\caption{The absolute error in the Gaussian approximations of $\varphi_{j}\left(x\right)$
for $j=1,2$ (plots (a) and (b)) , using $h=1$ and $M=200$.\label{fig:error for Gausssians, phi funcs}}

(a)~~~~\includegraphics[scale=0.5]{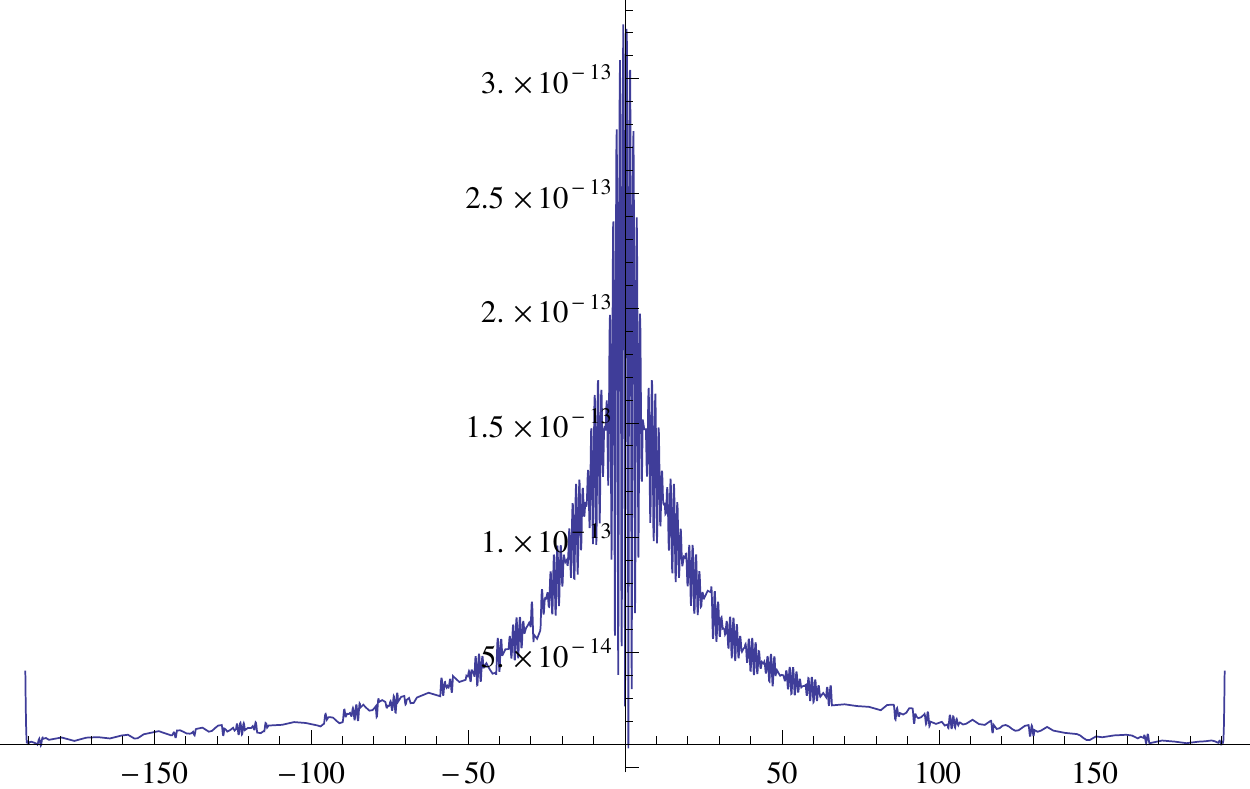}~~~~(a)~~~~\includegraphics[scale=0.5]{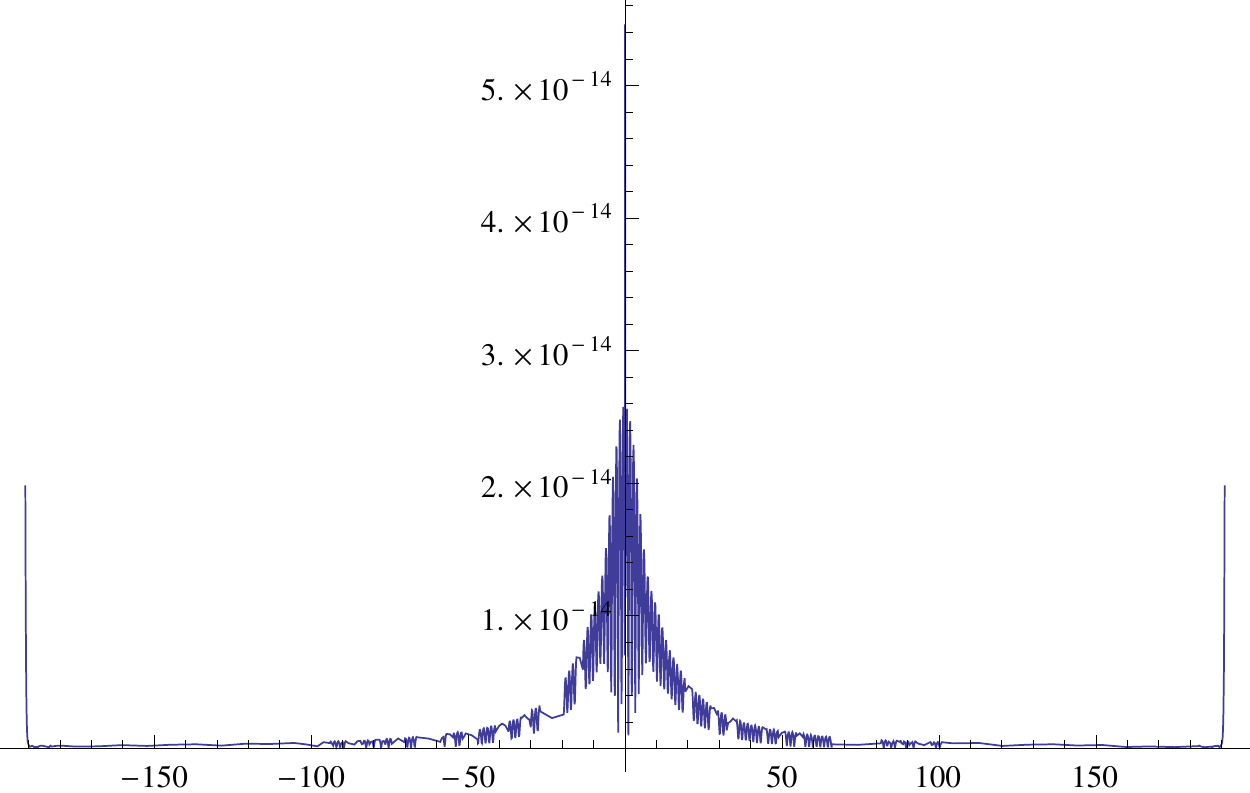}%
\end{figure}

In Figure~\ref{fig:error for Gausssians, phi funcs}, we plot the
error,\[
\left|\varphi_{j}\left(x\right)-\sum_{-\infty}^{\infty}c_{m,j}\psi\left(x+mh\right)\right|,\]
for the phi functions $\varphi_{1}\left(x\right)$ and $\varphi_{1}\left(x\right)$,
where we choose $h=1$ and $M=200$; notice that the choice of $h$
corresponds to the bandlimit of $\varphi_{j}\left(x\right)$. As shown
in Figure~\ref{fig:error for Gausssians, phi funcs}, the error is
smaller than $\approx3\times10^{-13}$ for all $-191\leq x\leq191$,
and is shown to begin to rise at the ends of the intervals, which
are close to $Mh$. This behavior can be understood by noting that
\[
\left|\varphi_{j}\left(x\right)-\sum_{-M}^{M}c_{m,j}\psi_{1}\left(x+m\right)\right|\leq\sum_{\left|m\right|>M}\left|c_{m,j}\right|\psi_{1}\left(x+m\right),\]
where we used that the support of $\widehat{\varphi_{j}}$ is contained
in $\left[-1/2,1.2\right]$. Since the functions $\psi_{1}\left(x+m\right)$
for $m>M$ decay rapidly away from $x=-m$, the error from truncation
is negligible when $\left|x\right|\leq\left(M-m_{0}\right)$ and $m_{0}=\mathcal{O}\left(1\right)$.

We remark that, for the function $e^{ix}$, it can
be shown (see the Appendix) that the approximation for $e^{ix}$ satisfies

\begin{equation}
\left|e^{ix}-\sum_{m=-M}^{M}c_{m}\psi_{h}\left(x+mh\right)\right|\leq\frac{1}{\widehat{\psi_{h}}\left(1\right)}\left(\sum_{k\neq0}\widehat{\psi_{h}}\left(\frac{k}{h}\right)+\sum_{\left|m\right|>M}\psi_{h}\left(x+mh\right)\right),\label{eq:error bound for e^ix, partial-1}\end{equation}
where $c_{m}$ is defined in (\ref{eq:coeff. for exp, gauss}). We
see that the first sum is negligible for e.g. $h\lesssim1$, owing
to the tight frequency localization of $\psi_{h}$. Similarly, the
second sum is negligible when $\left|x\right|\leq\left(M-m_{0}\right)h$
and $m_{0}=\mathcal{O}\left(1\right)$, owing to the tight spatial
localization of $\psi$.

\subsection{Rational approximation to a Gaussian\label{sub:Rational-approximation-to}}

We now discuss how to construct the approximation (\ref{eq:intro, rational approx to Gaussian}).

To do so, we first use AAK theory (see \cite{D-B-H-M:2013}
for details) to construct a near optimal rational approximation, \[
\left|\frac{1}{\sqrt{4\pi}}e^{-x^{2}/4}-\text{Re}\left(\sum_{j=1}^{N}\frac{b_{j}}{ix+\alpha_{j}}\right)\right|\leq\delta.\]
For an accuracy of $\delta\approx10^{-13}$, $13$ poles $\gamma_{j}$
are required.

\begin{figure}
\caption{\label{fig:Error-in-the-rational-approx-to-Gaussian}Error in the
rational approximation (\ref{eq:intro, rational approx to Gaussian})
to $e^{-x^{2}/4}$}

\includegraphics[scale=0.6]{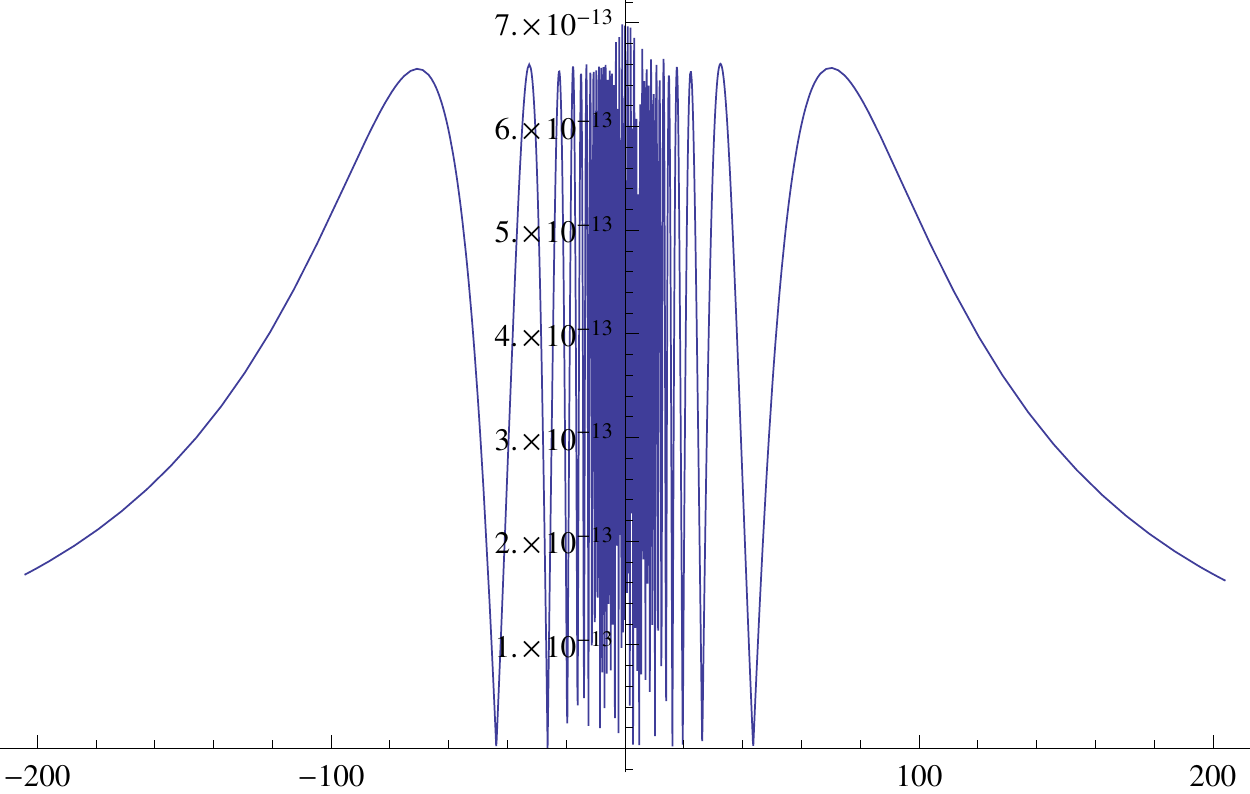}%
\end{figure}

Setting $\mu=\min_{j}\text{Re}\left(\alpha_{j}\right)$ , we next
look for a rational approximation to $\left(4\pi\right)^{-1/2}e^{-x^{2}/4}$
of the form 
\begin{equation} \label{eq:approx of Gaussian via rationals}
R\left(x\right)=\text{Re}\left(\sum_{j=-L}^{L}\frac{a_{j}}{ix+\mu+ij}\right),
\end{equation}
where we take $L=11$. We find the coefficients $a_{j}$ by minimizing
the $L^{\infty}$ error\[
\left\Vert \frac{1}{\sqrt{4\pi}}e^{-x^{2}/4}-\text{Re}\left(\sum_{j=-L}^{L}\frac{a_{j}}{ix_{n}+\mu+ij}\right)\right\Vert _{\infty},\]
where the points $x_{n}\in\left[-30,30\right]$ are chosen to be more
sparsely distributed outside the numerical support of $e^{-x^{2}/4}$;
the interval $\left[-30,30\right]$ is found experimentally to yield
high accuracy for the approximation over the entire real line. Finding
the coefficients $a_{j}$, $j=-L,\ldots,L$, that minimize the $L^{\infty}$
error can be cast as a convex optimization problem, and a standard
algorithm can be used (we use Mathematica). The resulting approximation
error is shown in Figure~\ref{fig:Error-in-the-rational-approx-to-Gaussian};
the error remains less than $\approx7\times10^{-13}$ for all $x\in\mathbb{R}$.

We display the real number $\mu$, and the coefficients $a_{j}$,
$j=1,\ldots,11$. In particular, these numbers are the only parameters
that are needed in order to construct rational approximations to general
functions on spatial intervals of any size.

\begin{table}
\caption{\label{tab:Coefficients for rational approx Gaussian}Coefficients
$a_{j}$, $j=-11,\ldots,11$, and number $\mu$, in the rational
approximation (\ref{eq:approx of Gaussian via rationals}).}

$\mu=-4.315321510875024$,

$a_{-11}=\left(-1.0845749544592896\times10^{-7},2.77075431662228\times10^{-8}\right)$,

$a_{-10}=\left(1.858753344202957\times10^{-8},-9.105375434750162\times10^{-7}\right)$,

$a_{-9}=\left(3.6743713227243024\times10^{-6},7.073284346322969\times10^{-7}\right)$,

$a_{-8}=\left(-2.7990058083347696\times10^{-6},0.0000112564827639346\right)$,

$a_{-7}=\left(0.000014918577548849352,-0.0000316278486761932\right)$,

$a_{-6}=\left(-0.0010751767283285608,-0.00047282220513073084\right)$,

$a_{-5}=\left(0.003816465653840016,0.017839810396560574\right)$,

$a_{-4}=\left(0.12124105653274578,-0.12327042473830248\right)$,

$a_{-3}=\left(-0.9774980792734348,-0.1877130220537587\right)$,

$a_{-2}=\left(1.3432866123333178,3.2034715228495942\right)$,

$a_{-1}=\left(4.072408546157305,-6.123755543580666\right)$,

$a_{0}=-9.442699917778205$,

$a_{1}=\left(4.072408620272648,6.123755841848161\right)$,

$a_{2}=\left(1.3432860877712938,-3.2034712658530275\right)$,

$a_{3}=\left(-0.9774985292598916,0.18771238018072134\right)$,

$a_{4}=\left(0.1212417070363373,0.12326987628935386\right)$,

$a_{5}=\left(0.0038169724770333343,-0.017839242222443888\right)$,

$a_{6}=\left(-0.0010756025812659208,0.0004731874917343858\right)$,

$a_{7}=\left(0.000014713754789095218,0.000031358475831136815\right)$,

$a_{8}=\left(-2.659323898804944\times10^{-6},-0.000011341571201752273\right)$,

$a_{9}=\left(3.6970377676364553\times10^{-6},-6.517457477594937\times10^{-7}\right)$,

$a_{10}=\left(3.883933649142257\times10^{-9},9.128496023863376\times10^{-7}\right)$,

$a_{11}=\left(-1.0816457995911385\times10^{-7},-2.954309729192276\times10^{-8}\right)$%
\end{table}

In Figure~\ref{fig:rational approx error, sin}, we show the resulting
rational approximations of $\cos\left(2\pi x\right)$ and $\sin\left(2\pi x\right)$,
which use the same $172$ complex-conjugate pairs of poles; the $L^{\infty}$
error is seen to be $\approx10^{-10}$ over the interval $-28\leq x\leq28$.

\begin{figure}
\caption{\label{fig:rational approx error, sin}Error in the rational approximations
of $\sin\left(2\pi x\right)$ and $\cos\left(2\pi x\right)$ (plots
(a) and (b)), for $-28\leq x\leq28$. These approximations use the
same $172$ pairs of complex-conjugate poles.}

(a)~~~~\includegraphics[scale=0.4]{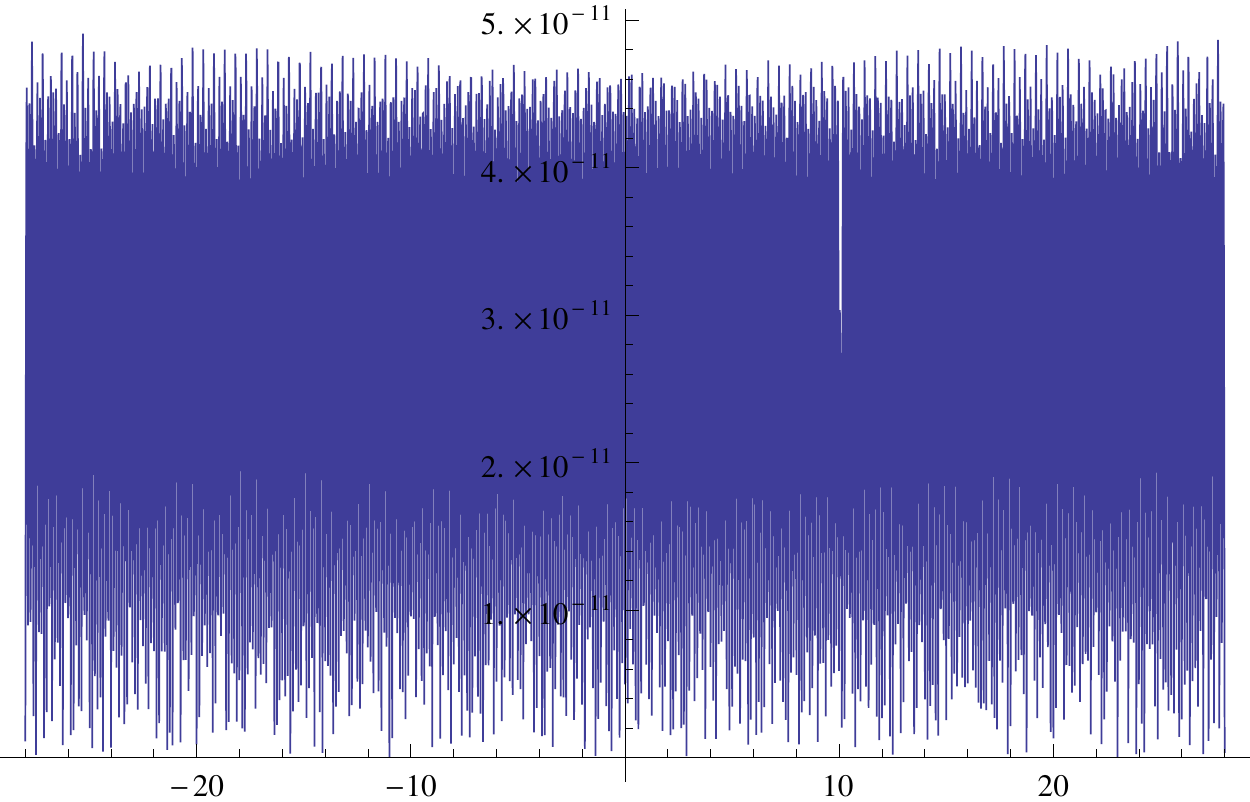}~~~~~(b)~~~~\includegraphics[scale=0.4]{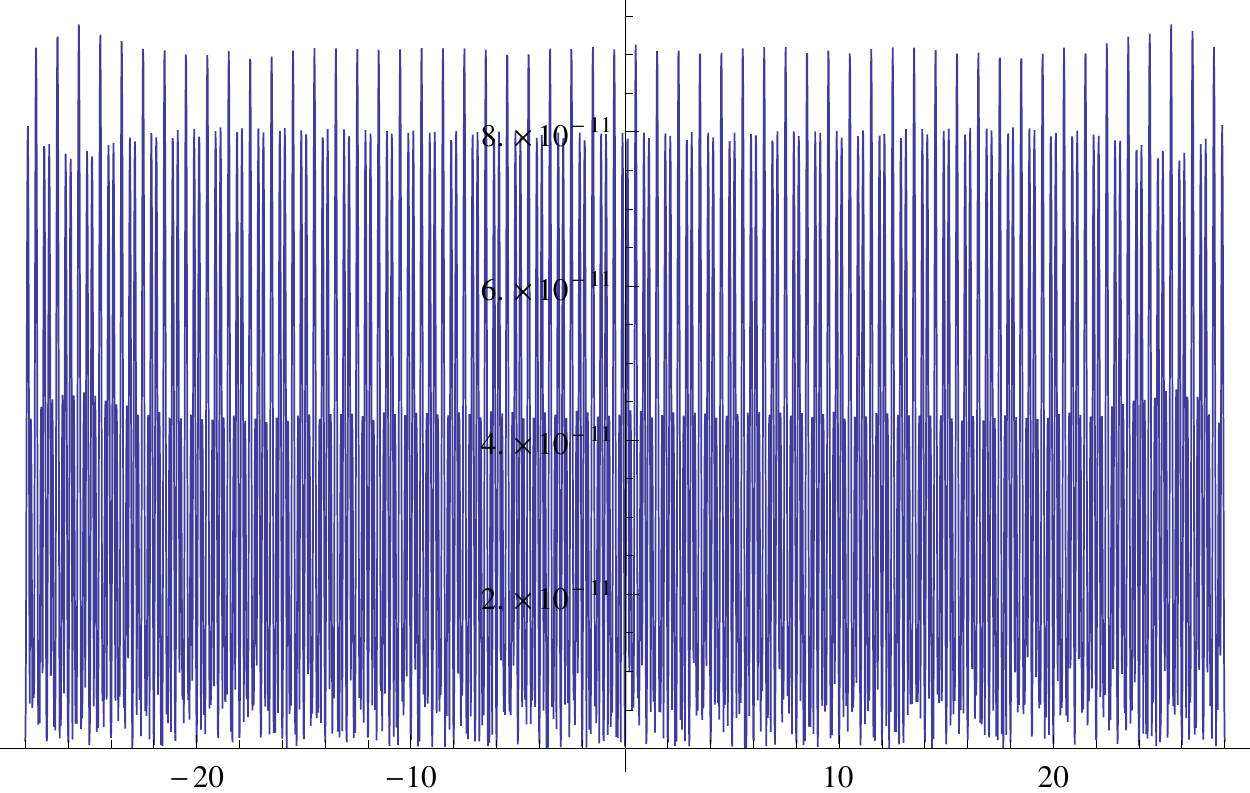}%
\end{figure}

\subsection{Constructing rational approximation of modulus bounded by unity} \label{Constructing rational approximation of modulus bounded by unity}

\begin{figure}
\caption{(a) Plot of the rational filter function $S\left(ix\right)$, for
$-60\leq x\leq60$.~~ (b) Plot of the difference $\left|S\left(ix\right)-1\right|$
for $-28\leq x\leq28$.\label{fig:partition of unity}}

~~~

(a)~~~~~~~~~~~~~~~~~~~~~~~~~~~~~~~~~~~~~~(b)

\includegraphics[scale=0.3]{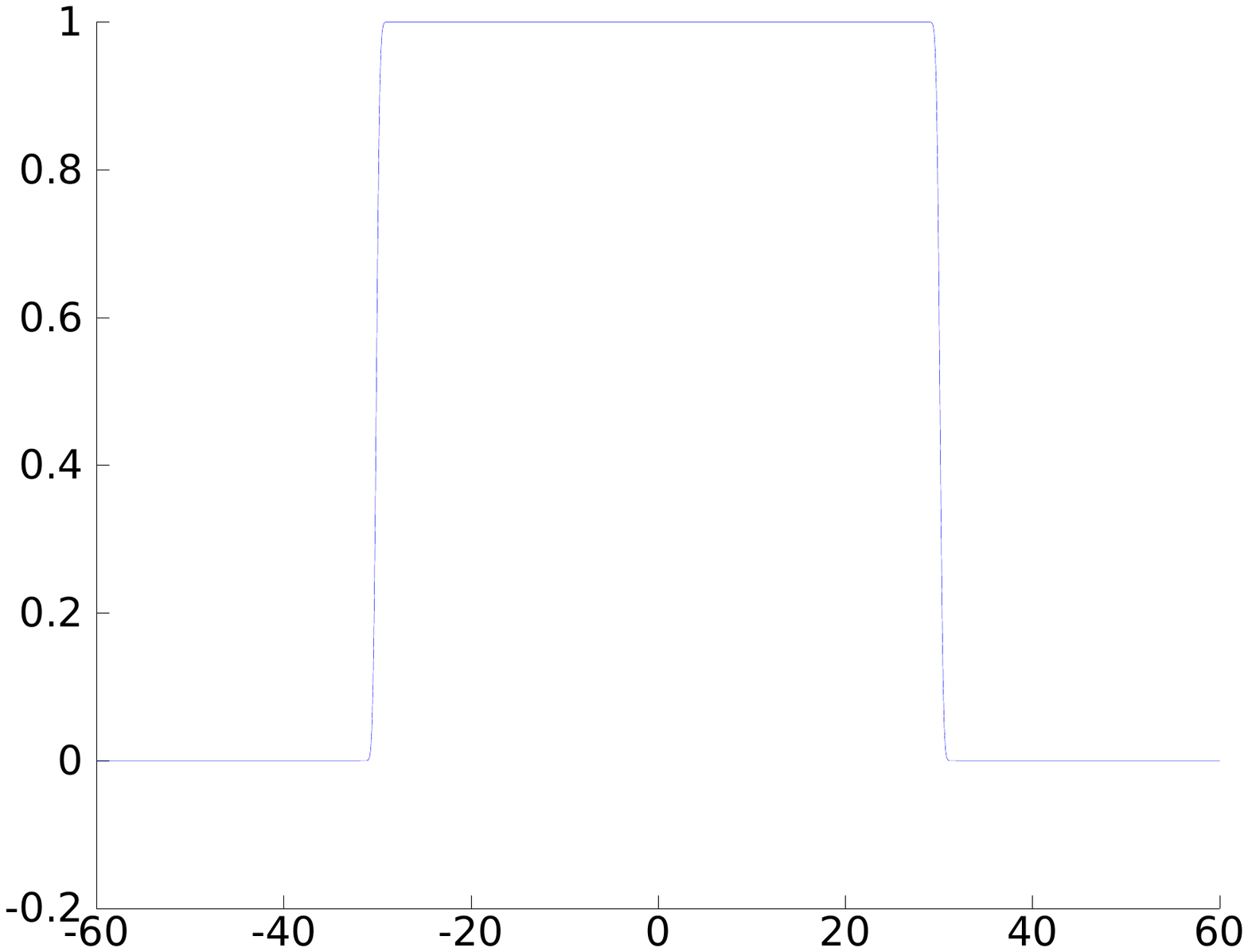}~~~\includegraphics[scale=0.3]{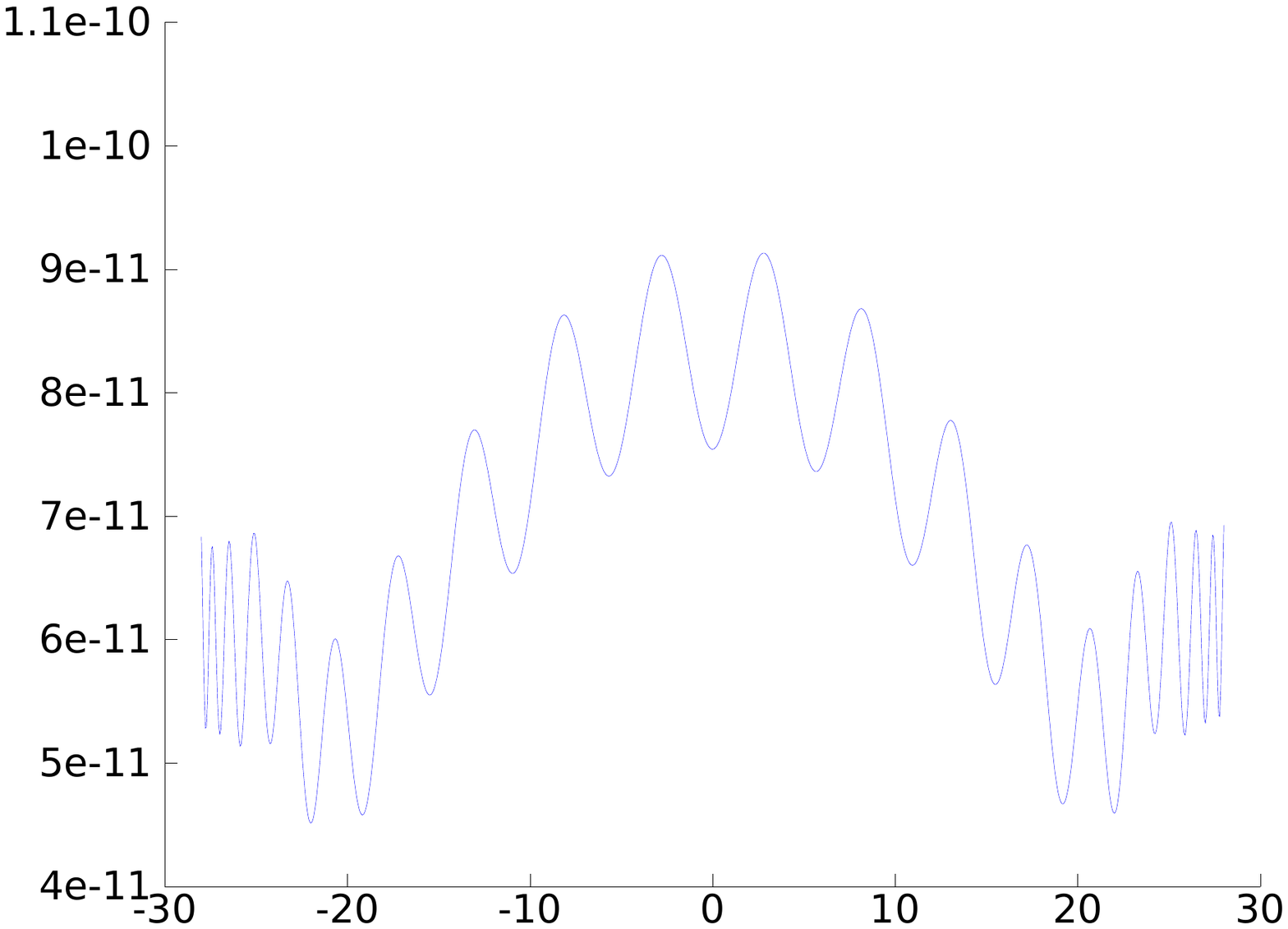}%
\end{figure}

For our applications, it is important that the approximation to $e^{ix}$
is bounded by unity on the real line. In particular, the Gaussian
approximation for $e^{ix}$ constructed in Section~\ref{sub:Gaussian-approximations-to}
has absolute value larger than one when $\left|x\right|\approx Mh$,
and this can lead to instability in repeated applications of $e^{t\mathcal{L}}$.

The basic idea is to construct a rational function $S\left(ix\right)$
that satisfies $S\left(ix\right)\approx1$ for $\left|x\right|\lesssim M_{0}h$
and $S\left(ix\right)\approx0$ for $\left|x\right|\gtrsim M_{0}h$.
As long as $M_{0}$ is slightly less than $M$, the function $S\left(ix\right)R_{M}\left(ix\right)$
accurately approximates $e^{ix}$ for $\left|x\right|\lesssim M_{0}h$,
and decays rapidly to zero for $\left|x\right|\gtrsim M_{0}h$. Therefore,
$\left|S\left(ix\right)R_{M}\left(ix\right)\right|\leq1$ for all
$x\in\mathbb{R}$, and repeated application of $S\left(t\mathcal{L}\right)R_{M}\left(t\mathcal{L}\right)\mathbf{u}_{0}$
is stable for all $t>0$. In Figure~\ref{fig:partition of unity},
we plot rational filter that uses $33$ complex-conjugate poles; we
see that $\left|S\left(ix\right)-1\right|\approx10^{-10}$ for $-28\leq x\leq28$.

Although the above approach results in a stable method, we have found
it more efficient to use a slightly modified version. This is motivated
by the following simple observation: since $\mathbf{u}_{0}\left(\mathbf{x}\right)$
is real-valued, \begin{equation}
\overline{\left(t\mathcal{L}-\alpha\right)^{-1}\mathbf{u}_{0}}=\left(t\mathcal{L}-\overline{\alpha}\right)^{-1}\mathbf{u}_{0}.\label{eq:conjugate, matrix inverses}\end{equation}
Recalling that the poles from Section~\ref{sub:Rational-approximation-to}
come in complex-conjugate pairs, only half the matrix inverses need
to be pre-computed and applied if (\ref{eq:conjugate, matrix inverses})
is used. However, directly using (\ref{eq:conjugate, matrix inverses})
results in numerical instabilities, where small errors in the high
frequencies are amplified after successive applications of $R_{M}\left(t\mathcal{L}\right)\mathbf{u}_{0}$.
The fix is to eliminate the errors in the high frequency components
by instead computing $S\left(k_{0}\Delta\right)R_{M}\left(t\mathcal{L}\right)\mathbf{u}_{0}$,
where $k_{0}$ is determined by the frequency content of $\mathbf{u}_{0}\left(\mathbf{x}\right)$
and the operator $S\left(k_{0}\Delta\right)$ only affects the highest
wavenumbers. Since the transition region between $S\left(ix\right)\approx1$
and $S\left(ix\right)\approx0$ can be made arbitrarily small (see
Figure~\ref{fig:partition of unity}), the operator $S\left(k_{0}\Delta\right)$
behaves like a spectral projector.

We now discuss how to construct $S\left(ix\right)$. To do so, we
use that (see \cite{MUL-VAR:2007}) \[
\left|\frac{1}{\widehat{\psi_{h}}\left(1\right)}\sum_{-\infty}^{\infty}\psi_{h}\left(x+hm\right)-1\right|\leq\frac{1}{h\widehat{\psi_{h}}\left(1\right)}\sum_{k\neq0}\widehat{\psi_{h}}\left(\frac{k}{h}\right),\]
which follows from the Poisson summation formula. For $h\lesssim1$,
the right hand side is negligible, owing to the tight frequency localization
of $\widehat{\psi_{h}}\left(\xi\right)$. Truncating the above sum
and using the tight spatial localiztion of $\psi_{h}\left(x\right)$,
we see that the function\begin{equation}
\chi\left(x\right)=\sum_{-M_{0}}^{M_{0}}\psi_{h}\left(x+mh\right),\label{eq:partition of unity-1}\end{equation}
is approximately unity for $\left|x\right|\lesssim M_{0}h$, and decays
to zero rapidly when $\left|x\right|\gtrsim M_{0}h$. It also holds
out that $\left|\chi\left(x\right)\right|\leq1$ for all $x\in\mathbb{R}$.
Therefore, using the techniques from Sections~\ref{sub:Gaussian-approximations-to}~and~\ref{sub:Rational-approximation-to},
we construct a rational approximation $Q\left(ix\right)$ to the function
$\chi\left(x\right)$ in (\ref{eq:partition of unity-1}), \begin{equation}
\left|Q\left(ix\right)-\sum_{-M_{0}}^{M_{0}}\psi_{h}\left(x+mh\right)\right|\leq\delta,\,\,\,\, x\in\mathbb{R},\label{eq:partition of unity}\end{equation}
 The number of poles required to represent the sub-optimal approximation
for $Q\left(x\right)$ can be drastically reduced with the reduction
algorithm \cite{HAU-BEY:2012}, which produces another proper rational
function $S\left(x\right)$ such that \[
\left|Q\left(ix\right)-S\left(ix\right)\right|\leq\delta_{0},\,\,\, x\in\mathbb{R},\]
and with a near optimally small number of poles for the prescribed
$L^{\infty}$ error $\delta_{0}$. Since the poles of $S\left(ix\right)$
and $R\left(ix\right)$ are distinct, the function $S\left(ix\right)R\left(ix\right)$
can be expressed as a proper rational function. The final function
$S\left(ix\right)$ is what is shown in Figure~\ref{fig:partition of unity}.

\section{Examples}
\label{sec:Examples}

\subsection{The 2D (rotating) shallow water equations}
\label{sub:Example-1}
We apply the technique proposed to the linear shallow water equations
\begin{eqnarray*}
\mathbf{v}_{t} & = & -fJ\mathbf{v}+\nabla\eta,\\
\eta_{t} & = & \nabla\cdot\mathbf{v},
\end{eqnarray*}
where all quantities are as in Section \ref{sec:reformulation}, cf.~equation (\ref{eq:RSW equations}).

\begin{figure}
\caption{\label{fig:errors for RSW, all} (a) Plots of the $L^{\infty}$ error,
$\left\Vert \mathbf{u}_{n}-e^{n\tau L}\mathbf{u}_{0}\right\Vert _{\infty}$,
versus the big time step $n\tau$, where $\tau=3$ and $1\leq n\le10$. Here the
approximation $\mathbf{u}_{n}$ is computed via
RK4, the Chebyshev polynomial method, and the rational approximation method. (b) Plots of
the computation time (min.) versus the big time step $n\tau$,  for the RK4, the Chebyshev polynomial method, and the rational approximation method.}
\
\

(a)

\includegraphics[scale=0.6]{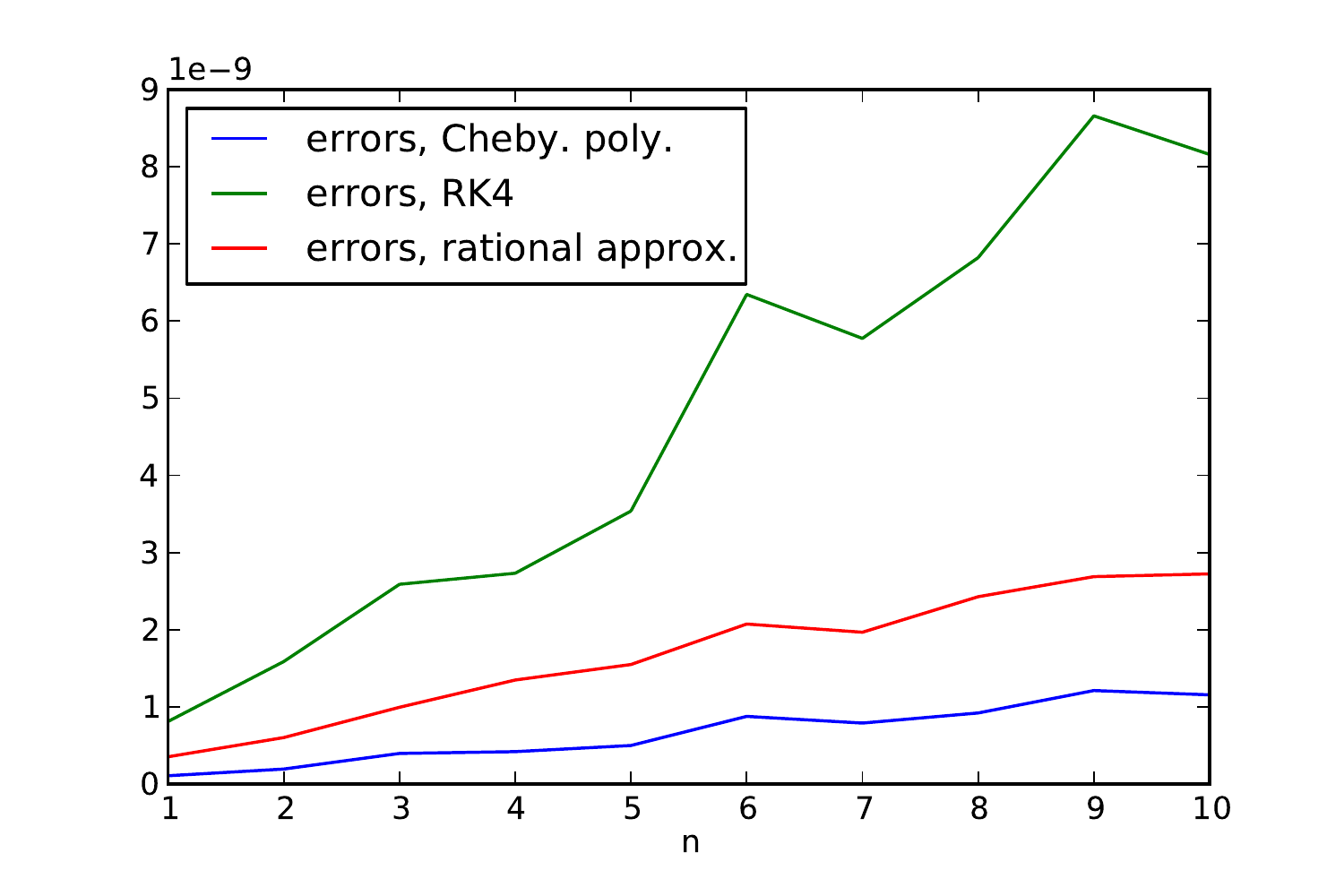}

(b)

\includegraphics[scale=0.6]{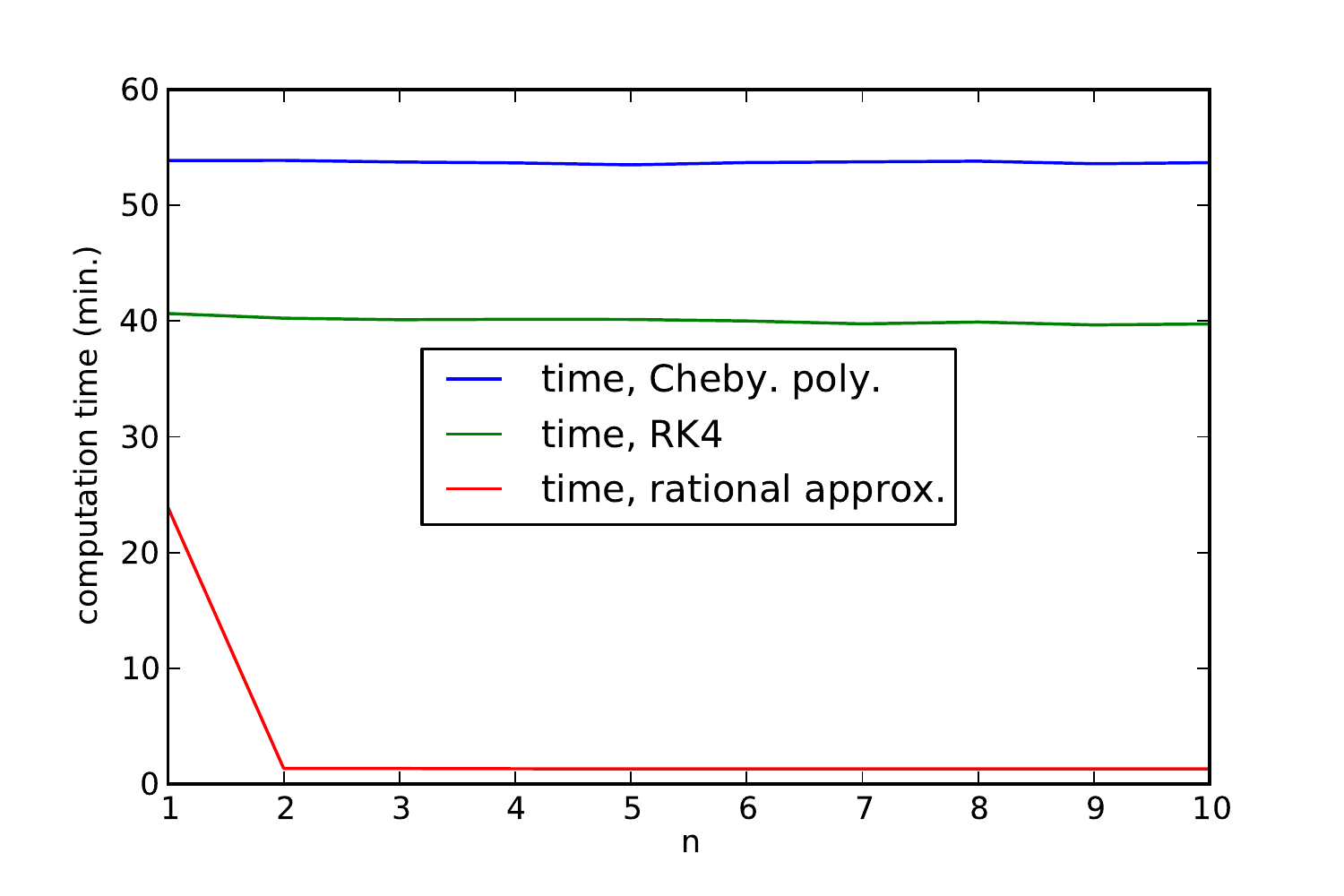}

\end{figure}

\begin{figure}
\caption{\label{fig:error for long time, RSW} Plot of the $L^{\infty}$ error,
$\left\Vert \mathbf{u}_{n}-e^{n\tau L}\mathbf{u}_{0}\right\Vert _{\infty}$,
versus the big time step $n\tau$, where $\tau=1$ and $1\leq n\le300$.
Here $\mathbf{u}_{n}$ denotes the numerical approximation to $e^{n\tau L}\mathbf{u}_{0}$,
as computed by the rational approximation (\ref{eq:error in exp(t*L)})
and the direct solver from Section~\ref{sec:Spectral-element-discretization}.}

\includegraphics[scale=0.5]{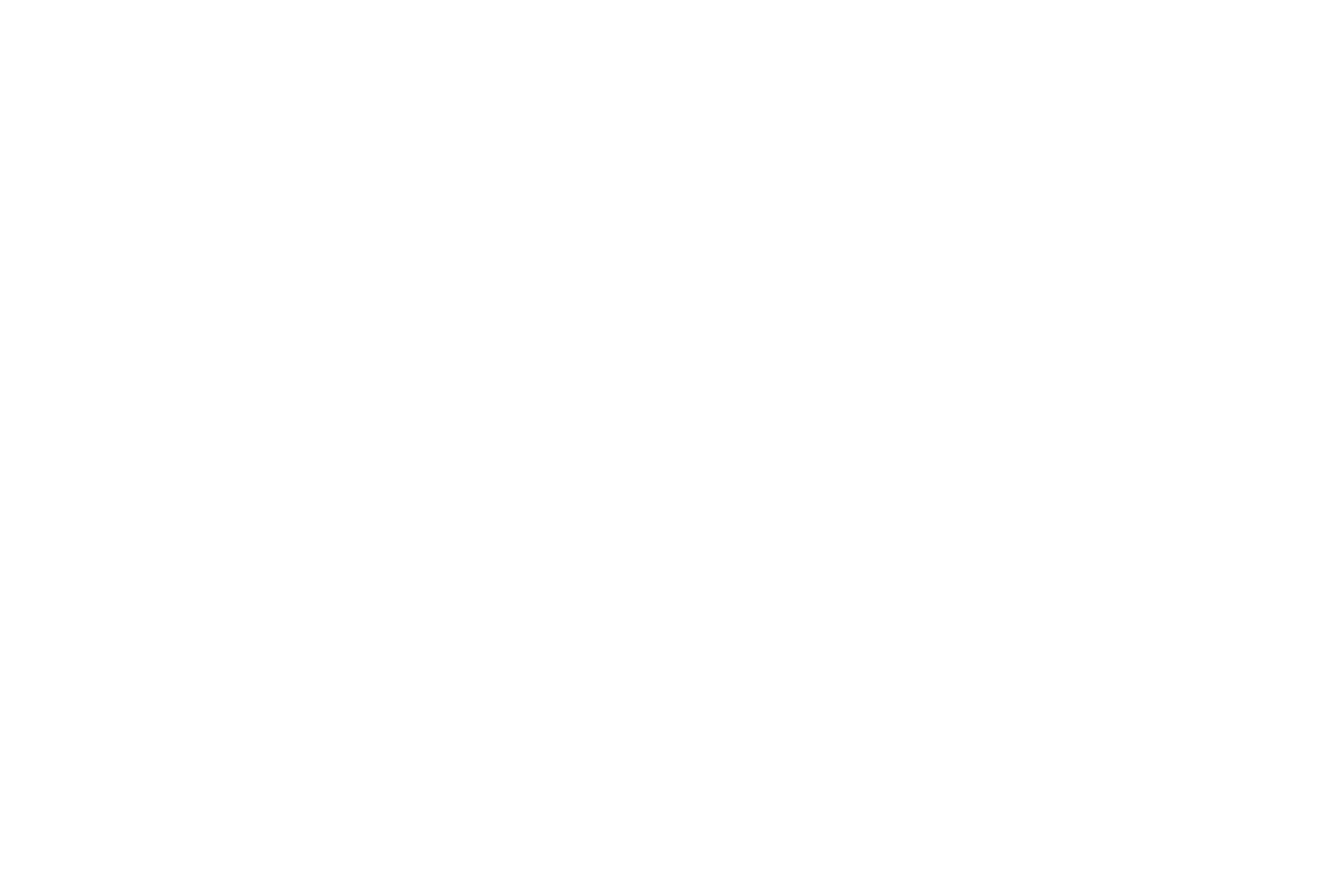}%
\end{figure}

\begin{table}
\label{tab:RSW, double res}
\caption{Comparison of the accuracy and efficiency of applying, $e^{\tau\mathcal{L}}\mathbf{u}_{0}$
and $\tau=1.5$, for system (\ref{eq:RSW equations}) and $\mathbf{u}_{0}$
in (\ref{eq:initial conditions, RSW eqns}). The comparison uses RK4,
Chebyshev polynomials, and the rational approximation (\ref{eq:error in exp(t*L)});
in the spatial discretization of all three comparisons, $12\times12=144$
elements and $16\times16=254$ Chebyshev quadrature nodes per element
are used. \label{tab:comparison for RSW}}

\begin{tabular}{|c|c|c|c|}
\hline
$e^{\tau\mathcal{L}}$, $\tau=1.5$ & $L^{\infty}$ error & time (min.) & pre-comp. (min.)\tabularnewline
\hline
Rational approx., & $2.1\times10^{-10}$ & $4.39$ & $103.1$\tabularnewline
$M=376$ terms &  &  & \tabularnewline
\hline
RK4 & $7.0\times10^{-10}$ & $131.9$ & NA\tabularnewline
 &  &  & \tabularnewline
\hline
Cheby. poly., & $1.1\times10^{-10}$ & $150.5$ & NA\tabularnewline
degree $12$  &  &  & \tabularnewline \hline
\end{tabular}%
\end{table}

We apply the algorithm in the spatial domain $\left[0,1\right]\times\left[0,1\right]$,
using periodic boundary conditions and a constant Coriolis force $f=1$ .
In this case, an exact solution can be computed analytically since the matrix
exponential is diagonalized in the Fourier domain, and can be rapidly applied
via the Fast Fourier Transform (FFT). In particular, \[
\mathcal{L}\left(\mathbf{r}_{\mathbf{k}}^{l}e^{i\mathbf{k}\cdot\mathbf{x}}\right)=i\omega_{\mathbf{k}}^{l}\mathbf{r}_{\mathbf{k}}^{l}e^{i\mathbf{k}\cdot\mathbf{x}},\]
where $\mathbf{r}_{\mathbf{k}}^{l}$ are eigenvectors of the matrix
\[
\left(\begin{array}{ccc}
0 & -f & igk_{1}\\
-f & 0 & igk_{2}\\
iHk_{1} & iHk_{2} & 0\end{array}\right),\]
and can be found in \cite{MAJDA:2003}.

We first compare the accuracy and efficiency of applying $e^{n\tau L}\mathbf{u}_{0}$,
for $\tau=3$ and $n=1,\ldots,10$, against $4$th order Runge-Kutta
(RK4) and against using Chebyshev polynomials. In particular, the
Chebyshev method uses the approximation\begin{equation}
e^{\Delta t\mathcal{L}}\mathbf{u}_{0}\approx J_{0}\left(i\right)\mathbf{u}_{0}+2\sum_{k=0}^{K}\left(i\right)^{k}J_{k}\left(-i\right)T_{k}\left(\Delta t\mathcal{L}\right)\mathbf{u}_{0},\label{eq:chebyshev expansion}\end{equation}
coupled with the standard recursion for applying $T_{k}\left(\Delta tL\right)$;
we choose a polynomial degree of $12$, which we find experimentally
is a good compromise between the time step size $\Delta t$ needed
for a given accuracy, and the number of applications of $\mathcal{L}$.
In all the time-stepping schemes, we use the same spectral element
discretization and parameter values as described above. All the algorithms
are implemented in Octave, including the direct solver described in
Section~\ref{sec:Spectral-element-discretization}.

\subsubsection{First test case for the shallow water equations}

We first consider the initial conditions\begin{eqnarray}
\eta\left(\mathbf{x}\right) & = & \sin\left(6\pi x\right)\cos\left(4\pi y\right)-\frac{1}{5}\cos\left(4\pi x\right)\sin\left(2\pi y\right),\nonumber \\
v_{1}\left(\mathbf{x}\right) & = & \cos\left(6\pi x\right)\cos\left(4\pi y\right)-4\sin\left(6\pi x\right)\sin\left(4\pi y\right),\nonumber \\
v_{2}\left(\mathbf{x}\right) & = & \cos\left(6\pi x\right)\cos\left(6\pi y\right).\label{eq:initial conditions, RSW eqns}\end{eqnarray}
For these initial conditions, we use $6\times6=36$ elements of equal
area, and $16\times16=256$ Chebyshev quadrature nodes for each element.
To assess the accuracy of the method, the exponential $e^{n\tau\mathcal{L}}\mathbf{u}_{0}$
is applied in the Fourier domain. When applying the operator exponential
using the rational approximation (\ref{eq:error in exp(t*L)}), we
use $M=376$ inverses and $\tau=3$ ; this results in an $L^{\infty}$
error of $3.4\times10^{-10}$ for a single (large) time step. For this choice of parameters in the spectral element discretization,
the cost of applying the solution operator of (\ref{eq:inverse oper, RSW eqns})---i.e.,
forming the right hand side of (\ref{eq:elliptic solve, RSW}), solving
(\ref{eq:elliptic solve, RSW}), and evaluating (\ref{eq:v_m, RSW})---is
about $4.5$ times more expensive than the cost of applying the forward
operator (\ref{eq:forward oper, RSW}) directly.


For the three time-stepping methods, the $L^{\infty}$ errors in the approximation of
$e^{n\tau\mathcal{L}}\mathbf{u}_{0}$, $n=1,\ldots,10$,
are plotted in Figure~\ref{fig:errors for RSW, all}, (a). Similarly, the total computation times (in minutes) of approximating
$e^{n\tau\mathcal{L}}\mathbf{u}_{0}$, $n=1,\ldots,10$,
are plotted in Figure~\ref{fig:errors for RSW, all}, (b)  (this includes the pre-computation time
for representing the inverses). From Figure~\ref{fig:errors for RSW, all}, (a), we see that the $L^{\infty}$ errors from
all three methods remain less than $10^{-8}$ for $n=1,\ldots,10$.  From Figure~\ref{fig:errors for RSW, all}, (b), we see that the first time step for
the rational approximation method is about half the cost of both RK4 and the Chebyshev polynomial method. However,
subsequent time steps for the new method is about $40$ times cheaper
than both RK4 and the Chebyshev polynomial method (for about the same accuracy).

\subsubsection{Second test case: doubling the spatial resolution} \label{subsec:Second test case: doubling the spatial resolution}

Next, we compute $e^{\tau\mathcal{L}}\mathbf{u}_{0}$, $\tau=1.5$,
with the initial conditions\begin{eqnarray}
\eta\left(\mathbf{x}\right) & = & \sin\left(12\pi x\right)\cos\left(8\pi y\right)-\frac{1}{5}\cos\left(8\pi x\right)\sin\left(4\pi y\right),\nonumber \\
v_{1}\left(\mathbf{x}\right) & = & \cos\left(12\pi x\right)\cos\left(8\pi y\right)-4\sin\left(12\pi x\right)\sin\left(8\pi y\right),\nonumber \\
v_{2}\left(\mathbf{x}\right) & = & \cos\left(12\pi x\right)\cos\left(12\pi y\right).\label{eq:initial conditions, RSW eqns-1}\end{eqnarray}
In particular, we double the bandlimit in each direction. In each
of the time-stepping schemes, we use $12\times12=144$ elements of
equal area, and $16\times16=256$ Chebyshev quadrature nodes for each
element. We again use $M=376$ inverses in (\ref{eq:error in exp(t*L)}).

We only examine the error and computation time for one big time step. For the rational approximation method,
we present both the pre-computation time for obtaining data-sparse representations of the $376$ inverses in (\ref{eq:error in exp(t*L)}),
and the computation time for applying the approximation in (\ref{eq:error in exp(t*L)}) (once the data-sparse representations are known).
The results are summarized in Table~\ref{tab:RSW, double res}. Since we only consider
a single time step, the pre-computation time and application time
are included separately. The main conclusion to draw from these results is
that doubling the spatial resolution does not appreciably change the relative efficiency
of the three time-stepping methods (once representations
for the inverse operators in (\ref{eq:error in exp(t*L)}) are pre-computed).

\subsubsection{Third test case: applying the operator exponential over a long time
interval}

Finally, we access the accuracy of the new method when repeatedly
applying $e^{\tau\mathcal{L}}$, $\tau=1$, in order to evolve the
solution over longer time intervals. In this example, we use the initial conditions

\begin{eqnarray}
\eta\left(\mathbf{x}\right) & = & \exp \left(-100\left( \left( x-1/2 \right)^2 +  \left( y-1/2 \right)^2  \right)  \right),\nonumber \\
v_{1}\left(\mathbf{x}\right) & = & \cos\left(6\pi x\right)\cos\left(4\pi y\right)-4\sin\left(6\pi x\right)\sin\left(4\pi y\right),\nonumber \\
v_{2}\left(\mathbf{x}\right) & = & \cos\left(6\pi x\right)\cos\left(6\pi y\right).\label{eq:initial conditions, RSW eqns, example 3}
\end{eqnarray}

Notice that these initial conditions cannot be expressed as a finite
sum of eigenfunctions of $\mathcal{L}$. We use the same spatial discretization parameters as in
Section~\ref{subsec:Second test case: doubling the spatial resolution}.

In Figure~\ref{fig:error for long time, RSW},
we show the $L^{\infty}$ error of the computed approximation $\mathbf{u}_{n}\left(\mathbf{x}\right)$
to $\mathbf{u}\left(\mathbf{x},n\tau\right)$, $n=1,\ldots,300$.
As expected, the error increases linearly in the number of applications
of the exponential. Notice that, due to the large step size of $\tau=1$,
the error accumulates slowly in time and the solution can be propagated
with high accuracy over a large number of characteristic wavelengths.

\subsection{Example 2}

In our second example, we consider the wave propagation problem \begin{equation}
u_{tt}=\kappa\Delta u,\,\,\,\,\mathbf{x}\in\left[0,1\right]\times\left[0,1\right],\label{eq:wave equation}\end{equation}
where $\kappa\left(\mathbf{x}\right)>0$ is a smooth function, the
initial conditions $u\left(\mathbf{x},0\right)$ and $u_{t}\left(\mathbf{x},0\right)$
are prescribed, and periodic boundary conditions are used.

\begin{table}
\caption{Comparison of the accuracy and efficiency for the operator exponential,
$e^{t\mathcal{L}}\mathbf{u}_{0}$ and $t=1.5$, for system (\ref{eq:forward oper, inhom})
and $\mathbf{u}_{0}$ in (\ref{eq:initial condition, wave prop}).
The comparison uses RK4, Chebyshev polynomials, and the rational approximation
(\ref{eq:error in exp(t*L)}); in the spatial discretization of all
three comparisons, $12\times12=144$ elements and $16\times16=254$
Chebyshev quadrature nodes per element are used. \label{tab:comparison for wave prop}}

\begin{tabular}{|c|c|c|c|}
\hline
$e^{t\mathcal{L}}$, $t=1.5$ & $L^{\infty}$ error & time (min.) & pre-comp. (min.)\tabularnewline
\hline
Rational approx., & $1.6 \times 10^{-9}$  & $3.76$ & $113.4$\tabularnewline
$M=376$ terms &  &  & \tabularnewline
\hline
RK4  & $3.5 \times 10^{-10}$ & $63.9$ & NA\tabularnewline
 &  &  & \tabularnewline
\hline
Cheby. poly., & $3.5 \times 10^{-8}$ & $57.5$ & NA\tabularnewline
degree $12$  &  &  & \tabularnewline \hline
\end{tabular}%
\end{table}

Since the procedure and results are
similar to those in Section~\ref{sub:Example-1}, we simply test the efficiency and accuracy of this
method over a single time step $\tau = 1.5$. In particular, we compare the accuracy and efficiency
for one application $e^{\tau\mathcal{L}}\mathbf{u}_{0}$, $\tau=1.5$,
against $4$th order Runge-Kutta (RK4) and against using Chebyshev
polynomials. In our numerical experiments, we use the initial condition\begin{equation}
u\left(x,y,0\right)=\sin\left(2\pi x\right)\sin\left(2\pi y\right)+\sin\left(4\pi x\right)\sin\left(4\pi y\right),\label{eq:initial condition, wave prop}\end{equation}
and $u_{t}\left(x,y,0\right)=0$. We also use \[
\kappa\left(x,y\right)=\left(\frac{3+\sin\left(4\pi x\right)}{4}\right)^{1/2}\left(\frac{3+\sin\left(4\pi y\right)}{4}\right)^{1/2}.\]
Finally, in the spatial discretization, we use $12\times12=144$ elements
with $16\times16=256$ points per element (for all three time-stepping
methods), and $M=376$ poles in (\ref{eq:error in exp(t*L)}). For
these parameters, the time to apply the inverse of (\ref{eq:operator for inhomog, inverse})---which
involves forming the right hand side in (\ref{eq:v for inhomog, inverse}),
solving for $v$, and computing (\ref{eq:w and z eqns})---is about
$5.2$ times more expensive than directly applying the forward operator
(\ref{eq:forward oper, inhom}).

Unlike Section~\ref{sub:Example-1}, the operator exponential is
not diagonalized in the Fourier domain. To assess the accuracy, we
use the Chebyshev polynomial method with a small enough step size
to yield an estimated error of less than $10^{-10}$. In particular,
we verify that the $L^{\infty}$ residual, $\left\Vert \mathbf{u}\left(\mathbf{x},t;\Delta t\right)-\mathbf{u}\left(\mathbf{x},t;\Delta t/2\right)\right\Vert _{\infty}$,
using numerical approximations to $\mathbf{u}\left(\mathbf{x},t\right)$  computed with step sizes $\Delta t$
and $\Delta t/2$ and the Chebyshev polynomial method, is less than
$10^{-10}$. We then use $\mathbf{u}\left(\mathbf{x},t;\Delta t/2\right)$
as a reference solution.

The results are summarized in Table~\ref{tab:comparison for wave prop}.
From this table, we see that the pre-computation time needed to represent
the $M=376$ solution operators in (\ref{eq:error in exp(t*L)}) is
$93$ minutes, and the computation time needed to apply the exponential
is $3.7$ minutes; the final accuracy in the $L^{\infty}$ norm is
given by $1.6 \times 10^{-9}$. For the Chebyshev polynomial method, $575$ time steps
of size $\Delta t \approx .0026$ are taken, for an overall time
of $57$ minutes; the final accuracy is given by $3.5 \times 10^{-8}$. Finally,
for RK4, $7,500$ time steps of size $\Delta t=1/5\times10^{-3}$
are taken, for an overall time of $63.9$ minutes; the final accuracy
is $3.5 \times 10^{-10}$.

\section{Generalizations}
\label{sec:generalizations}

The manuscript presents an efficient technique for explicitly computing a highly accurate approximation
to the operator $\varphi(\tau\mathcal{L})$ for the case where $\mathcal{L}$ is a skew-Hermitian operator
and where $\varphi(t) = e^{t}$, so that $\varphi(\tau\mathcal{L})$ is
the time-evolution operator of the hyperbolic PDE $\partial u /\partial t = \mathcal{L}u$. The technique
can be extended to more general functions $\varphi$. 
In particular, in using exponential integrators (cf. \cite{COX-MAT:2002}), 
it is desirable to apply functions  $\varphi_{j}\left(\tau \mathcal{L}\right)$, where $\varphi_{j} \left( \cdot \right)$ are the so-called phi-functions. In Section~\ref{sec:constructing-rational-approximations},
we presented (near) optimal rational approximations of the first few phi functions.
An important property of these representations is that the same poles can be used to simultaneously apply all the phi-functions, and with
a uniformly small error. In particular, linear combinations of the same $2 M +1$ solutions $\left(\tau \mathcal{L}-\alpha_{m}\right)^{-1}\mathbf{u}_{0}$,
$m=-M,\ldots,M$, can be used to apply  $\varphi_{j}\left(\tau \mathcal{L}\right)$ for $j=1,2,\ldots$.
In a similar way, linear combinations of the same $2 M +1$ solutions can be used to apply  $e^{s\mathcal{L}}$ for $0\leq s\leq \tau$.

In addition, where there is a priori knowledge of large spectral gaps---for example, when there is scale separation between fast and slow waves---the
techniques in this paper, coupled with those in \cite{HAU-BEY:2012}), can be used to construct
efficient rational approximations of $e^{i x}$ which are (approximately) nonzero only where the spectrum of $\mathcal{L}$ is nonzero. Since
suitably constructed rational approximations can capture functions with sharp transitions using
a small number of poles (see \cite{HAU-BEY:2012}), this approach requires
a potentially much smaller number of inverse applications.

\begin{appendix}

\section{Error bounds}  \label{Appendix: error bounds}

We now derive the error bound (\ref{eq:error bound for e^ix, partial-1}).
To do so, we use the Poisson summation formula,\[
\sum_{m=-\infty}^{\infty}\Psi_{h}\left(x+mh\right)=\frac{1}{h}\sum_{k=-\infty}^{\infty}e^{2\pi i\left(k/h\right)x}\widehat{\Psi_{h}}\left(\frac{k}{h}\right).\]
Applying this to $\Psi_{h}\left(x\right)=e^{-2\pi ix}\psi_{h}\left(x\right)$,
we have that\begin{eqnarray*}
\sum_{m=-\infty}^{\infty}\Psi_{h}\left(x+mh\right) & = & e^{-2\pi ix}\sum_{m=-\infty}^{\infty}e^{-2\pi imh}\psi_{h}\left(x+mh\right)\\
 & = & \frac{1}{h}\sum_{k=-\infty}^{\infty}e^{2\pi i\left(k/h\right)x}\widehat{\Psi_{h}}\left(\frac{k}{h}\right)\\
 & = & \frac{1}{h}\sum_{k=-\infty}^{\infty}e^{2\pi i\left(k/h\right)x}\widehat{\psi_{h}}\left(\frac{k}{h}+1\right),\end{eqnarray*}
where the last inequality uses the fact that \[
\widehat{\Psi_{h}}\left(\frac{k}{h}\right)=\widehat{\psi_{h}}\left(\frac{k}{h}+1\right).\]
%
{}Therefore,\[
\left|\sum_{m=-\infty}^{\infty}e^{-2\pi imh}\psi_{h}\left(x+mh\right)-\frac{\widehat{\psi_{h}}\left(1\right)}{h}e^{2\pi ix}\right|\leq\frac{1}{h}\sum_{k\neq0}\widehat{\psi_{h}}\left(\frac{k}{h}\right).\]
Finally, truncating the sum we obtain the bound (\ref{eq:error bound for e^ix, partial-1}).

{}

\end{appendix}

\bibliographystyle{plain}
\bibliography{common}

\end{document}